\documentclass[12pt,a4paper]{article}
\usepackage{latexsym,epsfig,fullpage,amsfonts,amssymb,amsmath,amscd,epic,graphics,rotating,theorem}
\usepackage[all]{xy}
\usepackage[usenames]{color}
\baselineskip 18pt \textwidth 15cm  \theoremstyle{plain}
\newtheorem{lemma}{Lemma}
\newtheorem{proposition}[lemma]{Proposition}
\newtheorem{remark}[lemma]{Remark}
\newtheorem{example}[lemma]{Example}
\newtheorem{theorem}{Theorem}

\newtheorem{corollary}{Corollary}
{\theorembodyfont{\rmfamily}
\begin{document}
\newcommand{\Pic}{{\operatorname{Pic}}}
\newcommand{\BH}{{\operatorname{bh}}}
\newcommand{\GW}{{\operatorname{GW}}}
\newcommand{\Conj}{{\operatorname{conj}}}
\newcommand{\Ev}{{\operatorname{Ev}}}
\newcommand{\idim}{{\operatorname{idim}}}
\newcommand{\DP}{{\operatorname{DP}}}
\newcommand{\Hom}{{\operatorname{Hom}}}
\newcommand{\Tors}{{\operatorname{Tors}}}
\newcommand{\Aut}{{\operatorname{Aut}}}
\newcommand{\Id}{{\operatorname{Id}}}
\newcommand{\Sing}{{\operatorname{Sing}}}
\newcommand{\codim}{{\operatorname{codim}}}
\newcommand{\Ann}{{\operatorname{Ann}}}
\newcommand{\ord}{{\operatorname{ord}}}
\newcommand{\mt}{{\operatorname{mt}}}
\newcommand{\Prec}{{\operatorname{Prec}}}
\newcommand{\proofend}{\hfill$\blacksquare$\bigskip}
\newcommand{\pperp}{{\perp{\hskip-0.4cm}\perp}}

\title{On higher genus Welschinger invariants of del Pezzo surfaces}


\author{Eugenii Shustin}

%
%
%
\date{}
\maketitle
\bigskip

\begin{abstract}
The Welschinger invariants of real rational algebraic surfaces
count real rational curves which represent a given divisor class and pass through a generic conjugation-invariant
configuration of points. No invariants counting real curves of positive genera are known in general.
We indicate particular situations, when Welschinger-type invariants
counting real curves of positive genera can be defined. We also prove the positivity and give
asymptotic estimates for such
Welschinger-type invariants for several del Pezzo surfaces of degree $\ge2$ and
suitable real nef and big divisor classes. In particular, this yields the existence of real curves
of given genus and of given divisor class
passing through any appropriate configuration of real points on the given surface.
\end{abstract}

\maketitle

\section{Introduction}

Welschinger invariants serve as genus zero open Gromov-Witten invariants. For
real rational symplectic manifolds \cite{Welschinger:2003,Welschinger:2005}, they count real rational pseudo-holomorphic curves,
realizing a given homology class, passing through a generic conjugation-invariant configuration of
points, and equipped with weights $\pm1$. In the case of real del Pezzo surfaces, Welschinger invariants
count real algebraic rational curves.
A more general approach used by J. Solomon
allowed him to define also invariants that count real curves of positive genera
with fixed complex and real structure of the normalization \cite[Theorem 1.3]{Solomon:2006}.
However, so far, no any general invariant way to count real curves of positive genera (without fixing their
complex and real structure) has been found.
In particular,
it follows from
\cite[Theorem 3.1]{Itenberg_Kharlamov_Shustin:2003}, that
if we do not fix complex and real structure of the normalization, then even the count of real plane elliptic curves
of any degree $\ge4$ equipped with Welschinger signs is {\it not invariant} of the choice of the point constraints.

The main goal of this note is to indicate situations, in which the ``bad" bifurcation of type
\cite[Theorem 3.1]{Itenberg_Kharlamov_Shustin:2003}
does not occur and in which Welschinger-type invariants
of positive genera can be defined. So, in Section \ref{pg}
we introduce higher genus invariants of
real del Pezzo surfaces with a disconnected real point set
and prove that
they indeed do not depend on the choice of point constraints and
on variation of the surface. In Section \ref{sec2}
we compute new invariants in several examples and exhibit a
series of invariants, which are positive and
are asymptotically comparable with Gromov-Witten
invariants. In particular, this yields the existence of real curves
of given genus and of given divisor class
passing through any appropriate configuration of real points on the given surface.

It is worth mentioning that
\cite[Theorem 1]{Itenberg_Kharlamov_Shustin:2009}
states that the count of tropical curves
of any genus with appropriate
tropical Welschinger signs is invariant of the choice of
tropical point constraints for any
toric surface. The reason why the ``bad" bifurcation
does not appear in the tropical limit is discussed
in \cite{Mikhalkin:2011,Mikhalkin:2011a}.

In our consideration we
intensively use techniques of \cite{Itenberg_Kharlamov_Shustin:2012}
and
\cite{Itenberg_Kharlamov_Shustin:2013}; for the
reader's convenience,
in Appendices 1 and 2, we present needed
statements from
these works in the form applicable to curves of arbitrary genus.

\section{Invariant count of real curves of positive
genera}\label{pg}

Let $X$ be a real del Pezzo surface with a nonempty
real point set ${\mathbb R} X$.
Denote by $\Pic^{\mathbb R}(X)\subset\Pic(X)$ the
subgroup of real divisor classes.
For any connected component $G\subset{\mathbb R} X$,
one can define a homomorphism $\BH_G:\Pic^{\mathbb R}(X)\to H_1(G;{\mathbb Z}/2)$
(cf. \cite{Borel_Haefliger:1961}),
which sends an
effective divisor class $D\in\Pic^{\mathbb R}(X)$ to the class
$[{\mathbb R} C\cap G]\in H_1(G;{\mathbb Z}/2)$, where $C\in|D|$
is any real curve. Indeed, it can
be viewed as the composition of the homomorphisms
$$H^{\Conj}_2(X)\to H_2(X/\Conj,{\mathbb R}X;{\mathbb Z}/2)\to
H_1({\mathbb R}X;{\mathbb Z}/2)
\to H_1(G;{\mathbb Z}/2)$$
given by $[\sigma]\mapsto[\sigma/\Conj]\mapsto[\partial(\sigma/
\Conj)]\mapsto[(\partial(\sigma/\Conj))
\cap G]$. It follows that,
for each $D\in\Pic^{\mathbb R}(X)$, there is a well-defined value
$\left(\BH_G(D)\right)^2\in\{0,1\}$.

Suppose that ${\mathbb R} X$ contains at least $g+1$
connected components $F_0,...,F_g$ for some $g\ge1$.
Put ${\widehat{F}}=F_0\cup...\cup F_g$, ${\underline{F}}=
(F_0,...,F_g)$. We say that a divisor class $D\in\Pic^{\mathbb R}(X)$ is
${\widehat{F}}$-compatible, if, for any connected component
$G\subset{\mathbb R} X\setminus{\widehat{F}}$, one has
$\BH_G(D)=0$. Note that ${\widehat{F}}$-compatible divisor
classes $D\in\Pic^{\mathbb R}(X)$ satisfy
$$DK_X\equiv\sum_{i=0}^g\left(\BH_{F_i}(D)\right)^2\mod2\ .$$
For any tuple $(r_0,...,r_g,m)$ of nonnegative integers,
introduce the space
${\mathcal P}_{{\underline{r}},m}(X,{\underline{F}})$
(where ${\underline{r}}=(r_0,...,r_g)$) of configurations
of $r_0+...+r_g+2m$ distinct points of $X$ such that $r_i$ of them
belong to $F_i$, $i=0,...,g$, and the others form $m$
complex conjugate pairs.

Choose any conjugation invariant class $\varphi\in
H_2(X\setminus{\widehat{F}};{\mathbb Z}/2)$ and
pick a big and nef, ${\widehat{F}}$-compatible
divisor class $D\in\Pic^{\mathbb R}(X)$ such that
\begin{equation}p_a(D)=(D^2+DK_X)/2+1\ge g\quad\text{and}\quad
-DK_X\ge g+1-\sum_{i=0}^g\left(\BH_{F_i}(D)\right)^2\ .
\label{e5}\end{equation} Then there exist nonnegative integers
$r_0,...,r_g,m$
such that
\begin{equation}
r_0+...+r_g+2m=-DK_X+g-1,\quad r_i\equiv\left(\BH_{F_i}(D)
\right)^2+1\mod2,\ i=0,...,g\ .
\label{e6}\end{equation} If $X$ is sufficiently generic
in its deformation class, and
${\boldsymbol w}\in{\mathcal P}_{{\underline{r}},m}(X,
{\underline{F}})$ is generic, then the set
${\mathcal C}^{\mathbb R}_g(X,D,{\boldsymbol w})$ of real
irreducible curve $C\in|D|$ of genus $g$, passing through
${\boldsymbol w}$, is finite and consists of only
immersed curves
(see Lemma \ref{ln12}).
Furthermore, each curve $C\in{\mathcal C}^{\mathbb R}_g
(X,D,{\boldsymbol w})$ has a one-dimensional real branch
in each of the components
$F_0,...,F_g$ of ${\mathbb R} X$. In particular, this yields
that $C\setminus{\mathbb R} C$ consists of two connected
complex conjugate components, and we denote one of them
by $C_{1/2}$. For any vector ${\underline{\varepsilon}}=
(\varepsilon_0,...,\varepsilon_g)$ with $\varepsilon_i=\pm1$, $i=0,...,g$, put
$$W_{g,{\underline{r}}}(X,D,{\underline{F}},
{\underline{\varepsilon}},\varphi,{\boldsymbol w})=
\sum_{C\in{\mathcal C}^{\mathbb R}_g(X,D,{\boldsymbol w})}(-1)^{s(C;{\underline{F}},{\underline{\varepsilon}})
+C_{1/2}\circ\varphi}\ ,$$
where $s(C,{\underline{F}},{\underline{\varepsilon}})$
is defined as follows:
if $C$ is nodal, then this is the number of those real
nodes of $C$ in ${\widehat{F}}$, which in $F_i$ are
represented in real local coordinates
as $x^2+\varepsilon_iy^2=0$, $i=0,...,g$ (a
node of type $x^2+y^2=0$ is called {\it solitary}, and of type
$x^2-y^2=0$ - {\it non-solitary}), and if $C$ is not nodal,
we locally deform each germ $(C,z)$, $z\in\Sing(C)$, moving its components
to a general position in an equivariant way, and then count real
nodes as in the nodal case. Since $C$ is immersed,
the number $s(C,{\underline{F}},{\underline{\varepsilon}})\mod2$
does not depend on the choice of local deformations of $C$.
Our main result is the following analog of Weslchinger's theorem
\cite{Welschinger:2003,Welschinger:2005} (see also \cite{Itenberg_Kharlamov_Shustin:2012}).

\begin{theorem}\label{t1} Let $X$ be a real del Pezzo surface,
$F_0,...,F_g$ connected components of ${\mathbb R} X$ for some $g\ge1$,
$D\in\Pic^{\mathbb R}(X)$ a nef and big,
${\widehat{F}}$-compatible divisor class satisfying
(\ref{e5}), $r_0,...,r_g,m$ nonnegative integers
satisfying (\ref{e6}), and $\varphi\in H_2(X\setminus
{\widehat{F}};{\mathbb Z}/2)$ a conjugation-invariant class.
Let ${\underline{\varepsilon}}=(\varepsilon_0,...,\varepsilon_g)$,
$\varepsilon_i=\pm1$, $i=0,...,g$. Then the following hold.

(1) The numbers $W_{g,{\underline{r}}}(X,D,{\underline{F}},
{\underline{\varepsilon}},\varphi,{\boldsymbol w})$ do not
depend on the choice of a generic
configuration ${\boldsymbol w}\in{\mathcal P}_{{\underline{r}},m}
(X,{\underline{F}})$ (which further on will be omitted in the notation).

(2) If tuples $(X,D,\underline{F},\varphi)$ and $(X',D',
\underline{F}',\varphi')$ are deformation equivalent,
then
$$
W_{g,{\underline{r}}}(X,D,{\underline{F}},{\underline{\varepsilon}},
\varphi)=W_{g,{\underline{r}}}(X',D',\underline{F}',{\underline{\varepsilon}},\varphi')\ .
$$
\end{theorem}

\begin{corollary}\label{c1}
Under the hypotheses of Theorem \ref{t1}, for any generic
configuration ${\boldsymbol w}\in{\mathcal P}_{{\underline{r}},m}
(X,{\underline{F}})$,
$$|W_{g,{\underline{r}}}(X,D,{\underline{F}},{\underline{\varepsilon}},
\varphi)|\le\#{\mathcal C}^{\mathbb R}_g(X,D,{\boldsymbol w})\le
 GW_g(X,D)\ ,$$ where $GW_g$ is the genus $g$ Gromov-Witten invariant.
In particular, if $W_{g,{\underline{r}}}(X,D,{\underline{F}},
{\underline{\varepsilon}},\varphi)\ne0$, then through any generic
configuration ${\boldsymbol w}\in{\mathcal P}_{{\underline{r}},m}(X,{\underline{F}})$, one
can trace a real curve $C\in|D|$ of genus $g$.
\end{corollary}

In Section \ref{sec2}, we exhibit several examples in which our invariants do not vanish, and therefore prove the existence of
real curves of positive genera passing through prescribed configurations of real points.

{\bf Proof of Theorem \ref{t1}.}
The proof follows the lines of \cite{Itenberg_Kharlamov_Shustin:2012}, where the case of rational curves
has been treated in full detail in the algebraic geometry framework. We only indicate principal
points of the argument, referring to Appendix 1, which contains all needed statements from
\cite{Itenberg_Kharlamov_Shustin:2012}.

The strategy of the proof is to verify that the studied enumerative numbers remain
constant in general variations
of the point constraints ${\boldsymbol w}\in{\mathcal P}_{{\underline{r}},m}(X,{\underline{F}})$ and of the surface $X$.

Let us fix a surface $X$ and consider
the space ${\mathcal P}^{\mathbb C}_n(X)$
of $n$-tuples of distinct points of $X$. Let $n=r_0+...+r_g+2m=-DK_X+g-1$. Then
${\mathcal P}_{{\underline{r}},m}(X,{\underline{F}})\subset{\mathcal P}^{\mathbb C}_n(X)$. Introduce the characteristic variety
$${\operatorname{Ch}}_n^{\mathbb C}(X,D)=\left\{{\boldsymbol w}\in{\mathcal P}_n^{\mathbb C}(X)
 \ \Bigg|\ \begin{array}{l}\text{there exists a Riemann surface}
\ S_g\ \text{of genus}\ g,\\
\text{an immersion}\ \nu:S_g\to X\ \text{and an}\ n\text{-tuple}\ {\boldsymbol p}\
\text{of distinct points of}\ S_g\\
\text{such that}\ \nu({\boldsymbol p})={\boldsymbol w},\ \nu(S_g)\in|D|,\ \text{and}\ h^1(S_g,{\mathcal N}_{S_g}^\nu(-
{\boldsymbol p}))>0\end{array}\right\} ,$$ where ${\mathcal N}_{S_g}^\nu=\nu^*{\mathcal T}X/{\mathcal T}S_g$ is the normal bundle.
If $p_a(D)>g$, this is a hypersurface in ${\mathcal P}^{\mathbb C}_n(X)$.
As pointed out in \cite[Theorem 3.1]{Itenberg_Kharlamov_Shustin:2003},
the invariance of
Welschinger numbers fails when the
(moving) configuration ${\boldsymbol w}$ hits ${\operatorname{Ch}}_n^{\mathbb C}(X,D)$.
Our key observation is that this event does not happen in our situation.

\begin{lemma}\label{l1}
Under the hypotheses of Theorem \ref{t1}, ${\mathcal P}_{{\underline{r}},m}(X,{\underline{F}})\cap{\operatorname{Ch}}_n^{\mathbb C}(X,D)=\emptyset$.
\end{lemma}

{\bf Proof.} Let $\nu:S_g\to X$ be a conjugation-invariant immersion, ${\boldsymbol p}\subset S_g$
be a conjugation-invariant
$n$-tuple such that $C=\nu^*S_g\in{\mathcal C}^{\mathbb R}_g(X,D,{\boldsymbol w})$, where
${\boldsymbol w}=\nu({\boldsymbol p})\in{\mathcal P}_{{\underline r},m}(X,{\underline F})$. Suppose that $h^1(S_g,{\mathcal N}_{S_g}^\nu(-{\boldsymbol p}))>0$. Then by Riemann-Roch
$h^0(S_g,{\mathcal N}_{S_g}^\nu(-{\boldsymbol p}))>0$. It is well known that $\nu_*{\mathcal N}_{S_g}^\nu={\mathcal J}_C^{cond}
\otimes{\mathcal O}_X(D)$, where ${\mathcal J}_C^{cond}
={\mathcal Ann}(\overline{\mathcal O}_C/{\mathcal O}_C)$
is the conductor ideal
sheaf on $C$ (see details, e.g., in
\cite[Section 4.2.4]{Dolgachev:2013}). Hence
$$H^0(C,{\mathcal J}_C^{cond}(-{\boldsymbol w})\otimes
{\mathcal O}_X(D))\simeq H^0(S_g,{\mathcal N}_{S_g}^\nu(-{\boldsymbol p}))\ne0$$ (in the case of
$w=z\in{\operatorname{Sing}}(C)$ for some
point $w\in{\boldsymbol w}$, we define the twisted sheaf
${\mathcal J}_C^{cond}(-{\boldsymbol w})$ as the limit when $w$ specializes to the point $z$ along a
component of the germ $(C,z)$). A real nonzero element of
$H^0(C,{\mathcal J}_C^{cond}(-{\boldsymbol w})\otimes{\mathcal O}_X(D))$ defines a real curve
$C'\ne C$ in the linear system $|D|$, which intersects $C$
at each singular point $z\in{\operatorname{Sing}}(C)$ with multiplicity
$\ge\delta(C,z)$ (see, \cite[Section 4.2.4]{Dolgachev:2013})
and in ${\boldsymbol w}$. In view of congruence
(\ref{e6}), $C'$ must intersect $C$ in (at least) one
additional point in each component $F_0,...,F_g$, and hence
\begin{eqnarray}CC'&\ge&\sum_{i=0}^g(r_i+1)+2m+2\delta(C)
\nonumber\\
&=&(-DK_X+g-1)+(g+1)+(D^2+DK_X+2-2g)=D^2+2\ ,\label{e2}
\end{eqnarray} which is a contradiction. \proofend

By Lemmas \ref{ln11}-\ref{ln16}, in a general smooth
equivariant deformation ${\boldsymbol w}_t$, $t\in[0,1]$,
of
${\boldsymbol w}={\boldsymbol w}_0$ in
${\mathcal P}_{{\underline{r}},m}(X,{\underline{F}})$, for $t$ in the complement to a finite set
$\Phi\subset[0,1]$, the curve collection ${\mathcal C}^{\mathbb R}_g(X,D,{\boldsymbol w})$
consists of immersed Riemann surfaces of genus $g$, and the values $t\in\Phi$
correspond to degeneration of some curves of ${\mathcal C}^{\mathbb R}_g(X,D,{\boldsymbol w})$
either into nonimmersed, birational images of Riemann surfaces of genus $g$, or into
curves listed in Lemma \ref{ln16}.
Lemmas \ref{l1} and \ref{p2X} yield that the numbers $W_{g,{\underline{r}}}(X,D,{\underline{F}},{\underline{\varepsilon}},\varphi,{\boldsymbol w})$
do not change as $t$ varies along any of the components of $[0,1]\setminus\Phi$. Next, we can suppose that
the configuration ${\boldsymbol w}_t$
is in general position on each of the finitely many curves $C=\nu(\hat C)$,
where $[\nu:\hat C\to X,{\boldsymbol p}]\in
{\mathcal C}^{\mathbb R}_g(X,D,{\boldsymbol w}_t)$, $t\in\Phi$. Then the constancy of the
numbers $W_{g,{\underline{r}}}(X,D,{\underline{F}},{\underline{\varepsilon}},\varphi,{\boldsymbol w})$
in these bifurcations follows from Lemmas \ref{p2X}, \ref{ln2}, and \ref{leg}.

The proof of statement (2) of Theorem \ref{t1} amounts in the verification of the constancy of
the number $W_{g,{\underline{r}}}(X,D,{\underline{F}},{\underline{\varepsilon}},
\varphi,{\boldsymbol w})$ when $X$ smoothly bifurcates through a uninodal
del Pezzo surface (see Section \ref{sec-dpudp} in Appendix 1)
The treatment is based on the use of an appropriate real version of the
Abramovich-Bertram-Vakil formula
\cite{Abramovich_Bertram:2001}, \cite[Theorem 4.2]{Vakil:2000}, and it literally
coincides with that in \cite[Section 4]{Itenberg_Kharlamov_Shustin:2013},
while the key points in this consideration are Lemmas \ref{ln12}(2ii) and
\ref{lem-abv}. \proofend

\section{Examples}\label{sec2}

\subsection{Small divisors}

\begin{proposition}\label{p1}
Suppose that the data $X,g,{\underline{F}},D,{\underline{r}},\varphi$
satisfy the hypotheses of Theorem \ref{t1}.

(1) If $p_a(D)=g$, then $W_{g,{\underline{r}}}(X,D,{\underline{F}},
{\underline{\varepsilon}},
\varphi)=(-1)^{C_{1/2}\circ\varphi}$, where $C$ is any smooth
curve from $|D|$.

(2) If $p_a(D)=g+1$, then
$$W_{g,{\underline{r}}}(X,D,{\underline{F}},{\underline{\varepsilon}},
\varphi)=\begin{cases}\sum_{i=0}^g\varepsilon_i(r_i+1-\chi(F_i)),\ &
\text{if}\ \varphi=0,\\
\sum_{i=0}^g\varepsilon_i(r_i+1-\chi(F_i))-\chi({\mathbb R} X\setminus
{\widehat{F}}),\ & \text{if}\ \varphi=[{\mathbb R} X\setminus{\widehat{F}}]
\end{cases}$$
\end{proposition}

{\bf Proof.} The first formula is evident, since the point constraints
define a unique smooth curve. In the second case,
the point constraints define a pencil of curves in $|D|$, which by
B\'ezout's argument similar to (\ref{e2}) have,
additionally to ${\boldsymbol w}$, an extra common point in each component
$F_0,...,F_g$, and hence the result follows from the Morse formula
after blowing up of all $\sum_{i=0}^g(r_i+1)$ real common points of the pencil.
\proofend

\begin{example}\label{e1}
Suppose that $X$ is a two-component real cubic surface in ${\mathbb P}^3$,
$F_0\simeq{\mathbb R} P^2$, $F_1\simeq S^2$, and let $g=1$.
Then (see \cite{Segre:1942}) $X$ contains precisely three real $(-1)$-curves
$E_1,E_2,E_3$ such that ${\mathbb R} E_1\cup{\mathbb R} E_2\cup{\mathbb R}
E_3\subset F_0$,
and each real affective, big and nef divisor can be represented as
$D=m_1E_1+m_2E_2+m_3E_3$ with $0<2m_i\le m_1+m_2+m_3$, $i=1,2,3$.
In particular, $-K_X=E_1+E_2+E_3$. Since $p_a(-2K_X-E_i)=2$, $i=1,2,3$, we have
$$W_{1,{\underline{r}}}(X,-2K_X-E_i,{\underline{F}},(\varepsilon_0,
\varepsilon_1),0)=\varepsilon_0r_0+\varepsilon_1(r_1-1)$$ for any
$r_0+r_1+2m=5$, $r_0\equiv0\mod2$, $r_1\equiv1\mod2$, $\varepsilon_0,
\varepsilon_1=\pm1$.
\end{example}

\subsection{Invariants of del Pezzo surfaces of degree $\ge2$}

Starting with the celebrated papers by Mikhalkin \cite{Mikhalkin:2005}
and Welschinger \cite{Welschinger:2005},
the problem of computation and analysis of the behavior of
(genus zero) Welschinger
invariants of real rational symplectic four-folds, in particular,
real del Pezzo surfaces
has been addressed in a series of papers (see, e.g.,
\cite{Arroyo_Brugalle_Medrano:2011,
Brugalle:2014,Brugalle_Puignau:2013,
Horev_Solomon:2012,Itenberg_Kharlamov_Shustin:2003,
Itenberg_Kharlamov_Shustin:2005,Itenberg_Kharlamov_Shustin:2009,
Itenberg_Kharlamov_Shustin:2013,Shustin:2006,
Welschinger:2007}). Some of the techniques developed there apply to
computation
of higher genus invariants introduced in Section \ref{pg}. In
this section, we demonstrate
examples of computations obtained by properly modified methods of
\cite{Itenberg_Kharlamov_Shustin:2013}. Similarly to
\cite{Itenberg_Kharlamov_Shustin:2013}, we stress on the positivity and asymptotic behavior of
our invariants, which particularly yield the existence of real curves of positive genus
passing through appropriate real point configurations.

Real del Pezzo surfaces are classified up to deformation equivalence by
their degree and the topology of the real point set (see
\cite{Degtyarev_Kharlamov:2002}).
In degree $\ge2$, we have the following surfaces $X$ with a
disconnected real part:
of degree $4$ with ${\mathbb R} X\simeq 2S^2$, of degree $3$ with
${\mathbb R} X\simeq{\mathbb R} P^2\pperp S^2$, and
of degree $2$ with ${\mathbb R} X\simeq{\mathbb R} P^2\pperp{\mathbb R} P^2$,
$({\mathbb R} P^2\#{\mathbb R} P^2)\pperp S^2$, $2S^2$, $3S^2$, or $4S^2$,
({\it cf.}, for instance, \cite[Section 5.1]{Itenberg_Kharlamov_Shustin:2013}).
For all of them, we can define elliptic invariants, for
the two last types invariants of genus $2$, and for the very
last one invariants of genus $3$.

\begin{proposition}\label{p2} Let $X$ be a real del Pezzo
surface of degree $\ge2$ such that ${\mathbb R} X$ contains (at least)
two connected components $F_0,F_1$
and let $D\in\Pic^{\mathbb R}(X)$ be a nef and big divisor class,
satisfying relations (\ref{e5}) for $g=1$. Then the following conditions are satisfied.
\begin{enumerate}\item[(i)] For any nonnegative integers $r_0,r_1$
satisfying (\ref{e6}) with $m=0$ and for any
conjugation-invariant class $\varphi\in H_2(X\setminus(F_0\cup F_1);
{\mathbb Z}/2)$, the invariants $W_{1,(r_0,r_1)}(X,D,(F_0,F_1),(1,1),
\varphi)$ do not depend on the choice
of the pair $r_0,r_1$
(thus, further on we omit subindex $(r_0,r_1)$ in the notation).
\item[(ii)] If $X$ is not of degree $2$ with ${\mathbb R} X\simeq2S^2$, then
\begin{equation}W_1(X,D,(F_0,F_1),(1,1),0)>0\ ,\label{e7}\end{equation} and
\begin{equation}\lim_{k\to\infty}\frac{\log W_1(X,kD,
(F_0,F_1),(1,1),0)}{k\log k}=\lim_{k\to\infty}
\frac{\log \GW_0(X,kD)}{k\log k}=-DK_X\ .\label{e10}\end{equation}
\item[(iii)] If $X$ is of degree $2$ with ${\mathbb R}
X\simeq2S^2$, then
\begin{equation}W_1(X,D,(F_0,F_1),(1,1),0)+W_{1,(-DK_X-1,1)}
(X,D,(F_0,F_1),(1,-1),0)>0\ ,\label{e12}
\end{equation} and
$$\lim_{k\to\infty}\frac{\log\Big(W_1(X,kD,(F_0,F_1),(1,1),0)+
W_{1,(-kDK_X-1,1)}(X,kD,(F_0,F_1),(1,-1),0)\Big)}
{k\log k}$$ \begin{equation}=\lim_{k\to\infty}\frac{\log
\GW_0(X,kD)}{k\log k}=-DK_X\ .\label{e13}
\end{equation}
\end{enumerate}
\end{proposition}

Statement (iii) of Proposition \ref{p2} can be generalized
to genus $2$ and $3$ invariants
of the surfaces $X$ of degree $2$ with ${\mathbb R} X\simeq3S^2$ or $4S^2$:

\begin{proposition}\label{p3}
(1) Let $X$ be a real del Pezzo surface of degree $2$ with
${\mathbb R} X\simeq 3S^2$ or $4S^2$, $F_0,F_1,F_2$ three distinct
connected components
of ${\mathbb R} X$,
$D\in\Pic^{\mathbb R}(X)$ a nef and big divisor class satisfying relation
(\ref{e5}) with $g=2$, $r_0,r_1$ odd positive integers satisfying
 $r_0+r_1=-DK_X$. Let ${\underline{r}}'=(r_0,r_1,1)$,
${\underline{F}}'=(F_0,F_1,F_2)$. Then
the invariant
$$W_{2,{\underline{r}}'}(X,D,{\underline{F}}',(1,1,\pm1)):=W_{2,
{\underline{r}}'}(X,D,{\underline{F}}',(1,1,1),0)+W_{2,
{\underline{r}}'}(X,D,{\underline{F}}',(1,1,-1),0)$$
does not depend on the choice of odd $r_0,r_1$ subject to $r_0+r_1=-DK_X$
(so, further on the subindex
${\underline{r}}'$ will be omitted), and it satisfies
$$W_2(X,D,{\underline{F}}',(1,1,\pm1))>0$$ and
$$\lim_{k\to\infty}\frac{\log W_2(X,kD,{\underline{F}}',(1,1,\pm1))}
{k\log k}=\lim_{k\to\infty}\frac{\log \GW_0(X,kD)}{k\log k}=-DK_X\ .$$

(2) Let $X$ be a real del Pezzo surface of degree $2$ with ${\mathbb R} X\simeq4S^2$, $F_0,F_1,F_2,F_3$
be `the connected components
of ${\mathbb R} X$, and
$D\in\Pic^{\mathbb R}(X)$ be a nef and big divisor class satisfying
relation (\ref{e5}) with $g=3$, $r_0,r_1$ odd positive integers satisfying
$r_0+r_1=-DK_X$.
Let ${\underline{r}}''=(r_1,r_2,1,1)$,
${\underline{F}}''=(F_0,F_1,F_2,F_3)$. Then
the invariant
$$W_{3,{\underline{r}}''}(X,D,{\underline{F}}'',(1,1,\pm1,\pm1)):=
\sum_{\varepsilon_2,\varepsilon_3=\pm1}W_{3,
{\underline{r}}''}(X,D,{\underline{F}}'',(1,1,\varepsilon_2,\varepsilon_3),0)$$
does not depend on the choice of odd $r_0,r_1$ subject to $r_0+r_1=-DK_X$ (so,
further on the subindex
${\underline{r}}''$ will be omitted), and it satisfies
$$W_3(X,D,{\underline{F}}'',(1,1,\pm1,\pm1))>0$$ and
$$\lim_{k\to\infty}\frac{\log W_3(X,kD,{\underline{F}}'',(1,1,\pm1,\pm1))}
{k\log k}=\lim_{k\to\infty}\frac{\log \GW_0(X,kD)}{k\log k}=-DK_X\ .$$
\end{proposition}

\begin{corollary}\label{c2}
(1) Under the hypotheses of Proposition \ref{p2}(ii) (respectively, \ref{p2}(iii)), through any
generic configuration ${\boldsymbol w}\in{\mathcal P}_{r_0,r_1,0}(X,(F_0,F_1))$
(respectively, ${\boldsymbol w}\in{\mathcal P}_{(-DK_X-1,1),0}(X,(F_0,F_1))$) one can draw
a real elliptic curve
$C\in|D|$ such that $C\supset{\boldsymbol w}$.

(2) Under the hypotheses of Proposition \ref{p3}(1) (respectively, \ref{p3}(2)), through any
generic configuration ${\boldsymbol w}\in{\mathcal P}_{{\underline{r}}',0}(X,
{\underline{F}}')$ (respectively, ${\mathcal P}_{{\underline{r}}'',0}(X,{\underline{F}}'')$)
one can draw a real curve
$C\in|D|$ of genus $2$ (respectively, $3$) such that $C\supset{\boldsymbol w}$.
\end{corollary}

\subsection{Proof of Proposition \ref{p2}}\label{sec2.3}

By blowing up suitable real points, we reduce the consideration to the only
surfaces of degree $2$. To treat this case, we use real versions of the
Abramovich-Bertram-Vakil formula and
the Caporaso-Harris-type formulas developed in \cite{Itenberg_Kharlamov_Shustin:2013},
as well as their direct extensions to elliptic curves.
We subsequently prove statements (i), (ii), and (iii).

\subsubsection{Proof of statement (i)}\label{sec-i}
Using Theorem
\ref{t1} and the construction of
\cite[Sections 4.2 and 5.2]{Itenberg_Kharlamov_Shustin:2013}, we can assume that
$X$ is a generic real fiber of an {\it elliptic ABV family} (in the terminology of
\cite[Section 5.2]{Itenberg_Kharlamov_Shustin:2013}), which is the following flat, conjugation-invariant family of surfaces
$\pi:{\mathfrak X}\to({\mathbb C},0)$: \begin{itemize}\item ${\mathfrak X}$ is a smooth three-fold;
\item
all fibers ${\mathfrak X}_t$, $t\ne0$, are del Pezzo of degree $2$; the fibers
${\mathfrak X}_t$, $t\in({\mathbb R},0)\setminus\{0\}$, are real, equivariantly deformation equivalent to $X$;
\item the central
fiber is ${\mathfrak X}_0=Y\cup Z$,
where $Y$ and $Z$ are smooth real surfaces transversally intersecting
along a smooth real rational curve $E$ which satisfies ${\mathbb R} E\ne\emptyset$ and
is such that $(Y,E)$ is a nodal
del Pezzo pair of degree $K_Y^2=2$ (i.e., $K_YE=0$, $-K_YC>0$ for any
irreducible curve $C\ne E$, and $(E^2)_Y=-2$, {\it cf.} \cite[Section 4]{Itenberg_Kharlamov_Shustin:2013}), $Z$ is a real quadric surface
with ${\mathbb R} Z\simeq S^2$, in
which $E$ is a hyperplane section (representing the divisor class $-K_Z/2$),
\item ${\mathbb R} E$ divides some connected component $F$ of ${\mathbb R} Y$ into
two parts $F_+,F_-$ so that the
components $(F_0)_t,(F_1)_t$ of ${\mathbb R}{\mathfrak X}_t$ (corresponding to the given components
$F_0,F_1$ of ${\mathbb R} X$), merge as $t\to0$ to $F_+$ and $F_-$, respectively.
\end{itemize} By \cite[Proposition 24]{Itenberg_Kharlamov_Shustin:2013},
$\Pic^{\mathbb R}(X)$ is naturally embedded into
$\Pic^{\mathbb R}(Y)$ as the orthogonal complement of $E$. Note also that the given class $\varphi\in
H_2(X\setminus(F_0\cup F_1);{\mathbb Z}/2)$ can be naturally identified with a conjugation-invariant
class in $H_2(Y\setminus F;{\mathbb Z}/2)$ (which we denote also by $\varphi$).

For a configuration ${\boldsymbol w}$ of $-DK_X=-DK_Y$ points in $F$ such that $r_0$
of them lie in $F_+$ and
$r_1$ other points lie in $F_-$ (we call such a configuration an {\it $(r_0,r_1)$-configuration}), denote by
${\mathcal C}^{\mathbb R}_1(Y,D,{\boldsymbol w})$ the set of real elliptic curves $C\in|D|_Y$ passing through ${\boldsymbol w}$.
By \cite[Proposition 2.1]{Shoval_Shustin:2013}, this is a finite set which consists of only immersed curves. Since
$DE=0$, any curve $C\in{\mathcal C}^{\mathbb R}_1(Y,D,{\boldsymbol w})$ has two one-dimensional real branches,
in particular, $C\setminus{\mathbb R} C$ splits into two connected components, one of which we denote by $C_{1/2}$.
Using \cite[Lemma 7]{Itenberg_Kharlamov_Shustin:2013}, we can replace each nonnodal singular point of any curve $C\in{\mathcal C}^{\mathbb R}_1(Y,D,{\boldsymbol w})$
by its local nodal equigeneric deformation and then correctly define the number
\begin{equation}
W_{1,(r_0,r_1)}(Y,D,F,\varphi,{\boldsymbol w})=\sum_{C\in{\mathcal C}^{\mathbb R}_1(Y,D,{\boldsymbol w})}(-1)^{s(C;F)+C_{1/2}\circ\varphi}\ ,\label{e4}
\end{equation}
where $s(C;F)$ is the number of solitary nodes of $C$ in $F$.

\begin{lemma}\label{l2} There exists a $(r_0,r_1)$-configuration ${\boldsymbol w}$ such that
$$W_{1,(r_0,r_1)}(X,D,(F_0,F_1),(1,1),\varphi)=W_{1,(r_0,r_1)}(Y,D,F,\varphi,{\boldsymbol w})\ .$$
\end{lemma}

{\bf Proof.} Take ${\boldsymbol w}$ to be a $D_0$-CH-configuration in the sense of Appendix 2, where $D_0\ge D$ is a suitable real effective divisor, and
$|{\boldsymbol w}\cap F_+|=r_0$,
$|{\boldsymbol w}\cap F_-|=r_1$.
Extend ${\boldsymbol w}$ up to a family of configurations
${\boldsymbol w}_t\subset{\mathbb R}{\mathfrak X}_t$, $t\in({\mathbb R},0)$, and note that
each elliptic curve $C_t\in{\mathcal C}_1^{\mathbb R}({\mathfrak X}_t,
D,{\boldsymbol w}_t)$ degenerates as $t\to 0$ either to
an elliptic curve $C_0\in{\mathcal C}_1^{\mathbb R}(Y,D,{\boldsymbol w})$, or
to the union of an elliptic curve
$C'_0\in{\mathcal C}_1^{\mathbb R}(Y,D-mE,{\boldsymbol w})$, $m>0$, and $2m$ generating
lines of the quadric $Z$ attached to $2m$ intersection points of $C'_0$ with $E$
({\it cf.} \cite[Lemma 22]{Itenberg_Kharlamov_Shustin:2013}). However, by
Lemma \ref{lem-ch}, all intersection points of $C'_0$ with $E$ are real, and hence
the above union of the generating lines of $Z$ is not real. Hence the latter degeneration of
$C_t$ is not possible, and we are done. \proofend

\begin{lemma}\label{l3}
If ${\boldsymbol w}$ is the $(r_0,r_1)$-configuration from Lemma \ref{l2} then
\begin{equation}
W_{1,(r_0,r_1)}(Y,D,F,\varphi,{\boldsymbol w})=W_{Y,E,\varphi+[{\mathbb R} Y\setminus F]}(D-E,0,2e_1,0)\ ,\label{e3}
\end{equation}
where the right-hand side is an ordinary $w$-number as defined in
\cite[Section 3.6]{Itenberg_Kharlamov_Shustin:2013}.
\end{lemma}

{\bf Proof.} By construction of Appendix 2, ${\boldsymbol w}=\{w_i\}_{i\in J}$, where
$J\subset\{1,...,N\}$, $|J|=r_0+r_1$. Let $k=\max J$. Consider degenerations of the curves
$C\in{\mathcal C}_1^{\mathbb R}(Y,D,{\boldsymbol w}$ induced by the deformation of
${\boldsymbol w}$, in which ${\boldsymbol w}'={\boldsymbol w}\setminus\{w_k\}$ stay fixed,
and $w_k$ specializes along the arc $L_k$ to the point $z_k\in E$ (see details in
Appendix 2).
By \cite[Proposition 2.6(2)]{Shoval_Shustin:2013}, any curve
$C\in{\mathcal C}^{\mathbb R}_1(Y,D,{\boldsymbol w})$ degenerates
\begin{enumerate}\item[(a)] either into the union $C'\cup E$, where $C'\in|D-E|$
is a real immersed elliptic curve,
passing through ${\boldsymbol w}'$,
intersecting $E$ at one point, and having there a smooth branch quadratically tangent to $E$,
\item[(b)] or into the union $C''\cup E$, where $C''\in|D-E|$ is a real immersed
rational curve, passing through
${\boldsymbol w}'$ and transversally intersecting $E$ in two distinct real points.
\end{enumerate} By \cite[Proposition 2.8(2)]{Shoval_Shustin:2013},
each curve $C'\cup E$ in item (a) gives rise to two
curves in ${\mathcal C}^{\mathbb R}_1(Y,D,{\boldsymbol w})$, which
are distinguished by (two) deformation patterns given
in\cite[Lemma 2.10(2)]{Shoval_Shustin:2013}, and which have opposite Welschinger signs (see \cite[Proposition 6.1(i)]{Shustin:2006}), and
therefore do not contribute to
$W_{1,(r_0,r_1)}(Y,D,F,\varphi,{\boldsymbol w})$. In its turn, each curve $C''\cup E$ in item (b) gives rise to
one curve in ${\mathcal C}^{\mathbb R}_1(Y,D,{\boldsymbol w})$. Furthermore, these curves $C''$ are counted by the number
$W_{Y,E,\varphi+[{\mathbb R} Y\setminus F]}(D,0,2e_1,0)$ with the same signs as the number
$W_{1,(r_0,r_1)}(Y,D,F,\varphi,{\boldsymbol w})$
counts the corresponding deformed curves in ${\mathcal C}^{\mathbb R}_1(Y,D,{\boldsymbol w})$ ({\it cf.} the right-hand sides in
(\ref{e4}) and \cite[Formulas (3) and (4)]{Itenberg_Kharlamov_Shustin:2013}), and hence (\ref{e3}) follows.
\proofend

Statement (i) of Proposition \ref{p2} is an immediate consequence of Lemmas \ref{l2} and \ref{l3}.

\begin{remark}\label{r1}
Lemmas \ref{l2} and \ref{l3} allow one to compute all considered invariants
$W_1(X,D,(F_0,F_1),(1,1),\varphi)$ via the recursive formula in
\cite[Theorem 2]{Itenberg_Kharlamov_Shustin:2013}. In Table \ref{tab1} we present
several values of this invariant for $D=-2K_X$ and various real del Pezzo surfaces $X$
(compared with the corresponding Gromov-Witten invariants of genus $1$,
{\it cf.} \cite[Examples
4.2 and 6.7]{Brugalle:2014}).

\begin{table}
\begin{center}
\begin{tabular}{|l||c|c|c|c|c|c|c|c|}
\hline $\deg X$ & 4 & 3 &
2& 2 & 2 & 2 & 2
  \\
\hline ${\mathbb R} X$ & {$2S^2$} & {${\mathbb R}P^2
\pperp S^2$} &
{$2{\mathbb R}P^2$}& {$({\mathbb R} P^2\#{\mathbb R} P^2)\pperp S^2$} & {$2S^2$}
& {$3S^2$} & {$4S^2$}
  \\
\hline
$W_1(X,-2K)$ & 112 & 36 & 12 & 12 & 4 & 8 & 16 \\
\hline
$GW_1(X,-2K)$ & 12300 & 1740 & 204 & 204 & 204 & 204 & 204 \\
\hline
\end{tabular}
\end{center}
\caption{Elliptic invariants of del Pezzo surfaces of degree $\ge2$}\label{tab1}
\end{table}

\end{remark}

\subsubsection{Proof of the positivity relation (\ref{e7})} By Lemmas \ref{l2} and \ref{l3}, to prove (\ref{e7}), it is enough to show that
\begin{equation}W_{Y,E,[{\mathbb R} Y\setminus F]}(D-E,0,2e_1,0)>0\ .\label{e9}\end{equation}
First, we prove an auxiliary inequality. Denote by ${\mathbb Z}_+^\infty$ the semigroup of vectors
$\alpha=(\alpha_1,\alpha_2,...)$ with countably many nonnegative integer coordinates such that
$\|\alpha\|=\sum_i\alpha_i<\infty$, and denote by ${\mathbb Z}_+^{\infty,odd}\subset{\mathbb Z}_+^\infty$ the subsemigroup
of vectors $\alpha$ such that $\alpha_{2i}=0$ for all $i\ge0$. By $e_1,e_2,...$ we be denote the standard
unit vectors in ${\mathbb Z}^\infty_+$. Put $I\alpha=\sum_{i\ge1}i\alpha_i$ for $\alpha\in{\mathbb Z}^\infty_+$.

\begin{lemma}\label{l4}
For any real nodal del Pezzo pair $(Y,E)$, introduced in Section
\ref{sec-i}, any nef divisor class
$D'\in\Pic^{\mathbb R}(Y)$ such that $D'E\ge0$ and $-D'K_Y>0$, and any vectors
$\alpha,\beta\in{\mathbb Z}_+^{\infty,odd}$ such that $I(\alpha+\beta)=D'E$, one has
\begin{equation}W_{Y,E,[{\mathbb R} Y\setminus F]}(D',\alpha,\beta,0)\ge0\ ,\label{e8}\end{equation}
where $W_{Y,E,[{\mathbb R} Y\setminus F]}(D',\alpha,\beta,0)$ is an ordinary $w$-number as defined
in \cite[Section 3.6]{Itenberg_Kharlamov_Shustin:2013}.
\end{lemma}

\begin{figure}
\begin{center}
\includegraphics[width=10cm]{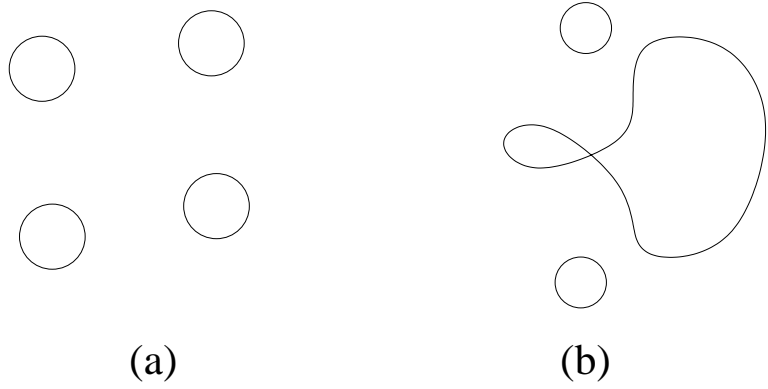}
\end{center}
\caption{Ramification quartics}\label{f1}
\end{figure}

{\bf Proof.}
For those pairs $(Y,E)$, which come from real del Pezzo surfaces $X$ with ${\mathbb R} X\simeq
S^2\pperp({\mathbb R} P^2\#{\mathbb R} P^2)$, ${\mathbb R} P^2\pperp{\mathbb R} P^2$, or $3S^2$, the claim follows from \cite[Lemma 39]{Itenberg_Kharlamov_Shustin:2013}.
Thus, we need to consider the only case of ${\mathbb R} X\simeq4S^2$. Via the anticanonical map $X\to{\mathbb P}^2$, the considered surface
$X$ is represented as the double covering of ${\mathbb P}^2$ ramified along a real smooth quartic curve
$Q_X$ having four ovals (see Figure \ref{f1}(a)), whereas ${\mathbb R} X$ doubly covers the four disks bounded by the ovals.
In turn, the family ${\mathcal X}$ can be obtained via the blow-up of the node of
the double covering of the trivial family ${\mathbb P}^2\times({\mathbb C},0)$ ramified along an
inscribed family of quartics with the nodal central quartic $Q_Y$ shown in Figure \ref{f1}(b).

To prove (\ref{e8}), we use induction on $R_Y(D',\beta):=-(K_Y+E)D'+\|\beta\|-1$. The base of induction is provided by
\cite[Proposition 9(1)]{Itenberg_Kharlamov_Shustin:2013}, where all nonzero values are equal to $1$. For the induction step, we apply the
suitably modified formula (6) from \cite[Theorem 2(2)]{Itenberg_Kharlamov_Shustin:2013}. In the left-hand side of
\cite[Formula (6)]{Itenberg_Kharlamov_Shustin:2013}, the summands of the first sum and the factors in the second sum, which
correspond to real divisor classes $D^{(i)}$ (in the notation of \cite{Itenberg_Kharlamov_Shustin:2013}), are nonnegative by the induction assumption,
whereas the factors corresponding to pairs of conjugate divisor classes may be negative. More precisely,
these factors correspond to pairs of conjugate $(-1)$-curves in $Y$ intersecting $E$. They can be viewed as follows
({\it cf.} \cite[Remark 23]{Itenberg_Kharlamov_Shustin:2013}):
there are exactly six
tangents to the quartic curve $Q_Y$ (Figure \ref{f1}(b)) passing through the node; they all are real,
and each one is covered by a pair of conjugate $(-1)$-curves in $Y$ intersecting in a real solitary node, which
projects to the tangency point
on $Q_Y$. Thus, a pair of $(-1)$-curves covering any of the two tangents to the real nodal branch of $Q_Y$ contributes
factor $(-1)$, while a pair of $(-1)$-curves covering any of the four tangents to the smooth ovals of $Q_Y$
contributes factor $1$.
Each summand of the second sum in the right-hand side of \cite[Formula (6)]{Itenberg_Kharlamov_Shustin:2013} can be
written as $(l+1)A_mB_{2l+m}$, where all the factors corresponding to pairs of
conjugate $(-1)$-curves are separated in $A_m$, where $m$ is the number of factors,
and the sum of the divisors classes appearing in the remaining part $B_{2l+m}$ equals
$D'-E-(2l+m)(K_Y+E)$. By \cite[Theorem 2(1g)]{Itenberg_Kharlamov_Shustin:2013}, any pair of $(-1)$-curves appears in
$A_m$ at most once. Thus, an easy computation converts \cite[Formula (6)]{Itenberg_Kharlamov_Shustin:2013} into
$$W_{Y,E,[{\mathbb R} Y\setminus F]}(D',\alpha,\beta,0)=\sum_{j\ge 1,\
\beta_j>0}W_{Y,E,[{\mathbb R} Y\setminus F]}(D',\alpha+e_j,\beta-e_j,0)+B_0+2B_1+B_2\ ,$$ which completes the proof in view of $B_0,B_1,B_2\ge0$ (by the induction assumption).
\proofend

Note that $D-E$ is nef (on $Y$). By \cite[Lemma 35(ii)]{Itenberg_Kharlamov_Shustin:2013}, it is enough to show that
$(D-E)E\ge0$ and $(D-E)E'\ge0$ for any $(-1)$-curve $E'$. We have $(D-E)E=DE-E^2=2$. For $(-1)$-curves
disjoint from $E$, we have $(D-E)E'=DE'\ge0$ by
the nefness of $D$. Any $(-1)$-curve $E'$ intersecting $E$ satisfies $E'E=1$,
and hence is not real (any real divisor has even intersection with $E$, since $[{\mathbb R} E]=0\in H_1({\mathbb R} Y)$).
Furthermore, $DE'>0$. Indeed, otherwise, $D$ would be disjoint both from $E'$ and from its complex conjugate
$\overline E'$; thus,
$D(E'+\overline E')=D(E+E'+\overline E')=0$, which in view of
$\max\{\dim|E'+\overline E'|,\dim|E+E'+\overline E'|\}=1$, would contradict the assumption $D^2>0$.
So, we conclude that $(D-E)E'=DE'-EE'=DE'-1\ge0$.

To complete the proof of (\ref{e7}), we establish a slightly stronger statement than (\ref{e9}).

\begin{lemma}\label{l5} For any real nodal del Pezzo pair $(Y,E)$ of degree
$\ge2$ with ${\mathbb R} E\ne\emptyset$ dividing
some connected component $F$ of ${\mathbb R} Y$, and any
nef divisor class $D'\in\Pic^{\mathbb R}(Y)$ such that $D'E=2$, one has
$$W_{Y,E,[{\mathbb R} Y\setminus F]}(D',0,2e_1,0)>0\ .$$
\end{lemma}

{\bf Proof.} We apply induction on $-D'K_Y$.

By \cite[Lemma 35(ii)]{Itenberg_Kharlamov_Shustin:2013}, $D'$ is nef on $X$. Since $D'\ne0$, it is effective on $X$,
and is presented by a smooth curve
(see, for instance, \cite[Theorems 3, 4, and Remark 3.1.4(B,C)]{Greuel_Lossen_Shustin:1998}, where the condition
$p_a(D)\ge0$ trivially follows from \cite[Formula (3.1.2)]{Greuel_Lossen_Shustin:1998}), and hence $-D'K_Y=-D'K_X>0$.
Furthermore, $-D'K_Y\ne1$. Indeed, otherwise, by the genus formula $(D')^2\equiv-D'K_Y=1\mod2$, that is $(D')^2\ge1$,
and thus, $p_a(D')\ge1$. However, $-D'K_X=1$ and $\dim|K_X|\ge1$ would imply that a general
curve $C\in|D'|_X$ is rational, which is a contradiction. Hence $-D'K_X\ge2$. Suppose that
$-D'K_X=-D'K_Y=2$. This yields $-D'(K_Y+E)=0$, which ({\it cf.} \cite[Lemma 35(iii)]{Itenberg_Kharlamov_Shustin:2013}) leaves the only case
$K_Y^2=2$ and $D'=-K_Y-E$, represented by a smooth rational curve, which finally yields
$W_{Y,E,[{\mathbb R} Y\setminus F]}(D',0,2e_1,0)=1$.

Suppose that $-D'K_Y>2$. By the genus formula, $(D')^2>0$. Then
$D'E'>0$ for any $(-1)$-curve $E'$ intersecting $E$ ({\it cf.} the argument
in the proof of the nefness of $D-E$ above). If $D'$ is disjoint from a real $(-1)$-curve $E'$ such that
$E'E=0$, we blow down $E'$. If $D'$ is disjoint from a nonreal $(-1)$-curve $E'$ such that $E'E=0$, then
$E'\overline E'=0$ (since otherwise $D'$ would be disjoint with curves in the one-dimensional linear system
$|E'+\overline E'|$ contrary to $(D')^2>0$), and then we blow down both $E'$ and $\overline E'$. After finitely many such steps
we arrive to a real nodal del Pezzo surface $(Y',E)$ of degree $\ge 2$ and a nef and big divisor
class $D'\in\Pic^{\mathbb R}(Y')$ such
that $D'E=2$, $-D'K_{Y'}=D'K_Y$, and $D'E'>0$ for any
$(-1)$-curve in $Y'$. It follows that $(D'+K_{Y'})E=2$, and that $D'+K_{Y'}$ nonnegatively intersects any $(-1)$-curve on
$Y'$. Hence $D'+K_{Y'}$ is nef on $Y'$. Since $-(D'+K_{Y'})K_{Y'}<-D'K_{Y'}=-DK_Y$, we have
$W_{Y',E,[{\mathbb R} Y'\setminus F']}(D'+K_{Y'},0,2e_1,0)>0$, where $F'\subset{\mathbb R} Y'$ is the image of $F$. Then, by
\cite[Formula (6)]{Itenberg_Kharlamov_Shustin:2013} and by Lemma \ref{l4},
$$W_{Y,E,[{\mathbb R} Y\setminus F]}(D',0,2e_1,0)=W_{Y',E,[{\mathbb R} Y'\setminus F']}(D',0,2e_1,0)$$ $$\ge
W_{Y',E,[{\mathbb R} Y'\setminus F']}(D'+K_{Y'},0,2e_1,0)\cdot W_{Y',E,[{\mathbb R} Y'\setminus F']}(-K_{Y'}-E,0,2e_1,0)>0\ ,$$
where $W_{Y',E,[{\mathbb R} Y'\setminus F']}(-K_{Y'}-E,0,2e_1,0)=1$, because $p_a(-K_{Y'}-E)=0$, and hence
a general curve in $|-K_{Y'}-E|_{Y'}$ is smooth rational. \proofend

\subsection{Proof of the asymptotic relation (\ref{e10})} It is enough to show that
\begin{equation}\log W_1(X,kD,(F_0,F_1),(1,1),0)\ge(-DK_X)k\log k+O(k)\ ,
\label{e11}\end{equation} since by Lemmas \ref{l2} and
\ref{l3}, and by \cite[Theorem 1]{Itenberg_Kharlamov_Shustin:2005},
$$\log W_1(X,kD,(F_0,F_1),(1,1),0)=\log W_{Y,E,[{\mathbb R} Y\setminus F]}(kD-E,0,2e_1,0)$$
$$\le \log\GW_0(X,kD)=(-DK_X)k\log k+O(k)\ .$$ Using Lemmas \ref{l4} and \ref{l5}, and
\cite[Formula (6)]{Itenberg_Kharlamov_Shustin:2013}, we derive for any $k\ge2$
$$W_{*}(kD-E,0,2e_1,0)\ge\frac{1}{2}\cdot\frac{(-kDK_Y-2)!}{(-iDK_Y-1)!(-(k-i)DK_Y-1)!}$$
$$\times\sum_{i=1}^{k-1}\Big[4\cdot
W_{*}(iD-E,0,2e_1,0)\cdot W_{*}((k-i)D-E,0,2e_1,0)\Big]\ .$$
where the asterisk stands for the subindex $(Y,E,[{\mathbb R} Y\setminus F])$. This inequality yields that the positive sequence
$$a_n=\frac{W_{*}(nD-E,0,2e_1,0)}{(-nDK_Y)!},\quad n\ge 1\ ,$$ satisfies the relation
$a_n\ge\lambda\sum_{i=1}^{n-1}a_i$ with some absolute constant $\lambda>0$. By
\cite[Lemma 38]{Itenberg_Kharlamov_Shustin:2013}, $a_n\ge\xi_1\xi_2^n$, $n\ge1$, with some positive $\xi_1,\xi_2$, which leads to
(\ref{e11}).

\subsubsection{Proof of statement (iii)}\label{sec-iii} Let $X$ be a real del Pezzo surface of degree $2$ with ${\mathbb R} X\simeq 2S^2$ and
$D\in\Pic^{\mathbb R}(X)$ a nef and big divisor class. So, $F_0\simeq F_1\simeq S^2$,
and we let $r_0=-DK_X-1$, $r_1=1$. Since all such surfaces are equivariantly deformation equivalent
and in view of Theorem \ref{t1}(2),
we can suppose that $X$ is a fiber ${\mathfrak X}'_\tau$, $\tau>0$, if a flat
conjugation-invariant
family ${\mathfrak X}'\to({\mathbb C},0)$ of surfaces, along which the component $F_1$
collapses to an isolated real nodal point so that in a neighborhood of the node
the family is representable as $x_1^2+x_2^2+x_3^2=\tau$. Following \cite[Section 4.2]{Itenberg_Kharlamov_Shustin:2013},
we perform the base change $\tau=t^2$ and blow up the node obtaining finally a conjugation-invariant
family (called a {\it $3$-unscrew} ${\mathfrak X}\to({\mathbb C},0)$ in
\cite[Section 4.2]{Itenberg_Kharlamov_Shustin:2013})
with the central fiber ${\mathfrak X}_0=Y\cup Z$, where $E=Y\cap Z$ is a smooth real
rational curve with ${\mathbb R} E=\emptyset$, $(Y,E)$ being a real nodal del Pezzo pair with ${\mathbb R} Y\simeq S^2$, and
$Z$ is a quadric surface in which $E$ represents the divisor class $-K_Z/2$ and which
has the real part ${\mathbb R} Z\simeq S^2$. Pick a generic point $w_0\in{\mathbb R} Z$ and a generic configuration ${\boldsymbol w}'\subset{\mathbb R} Y$
of $-DK_Y-1=-DK_X-1$ distinct points in ${\mathbb R} Y$, and extend $\{w_0\}\cup{\boldsymbol w}'$ to smooth equivariant sections
$t\mapsto{\boldsymbol w}_t$ of
the family ${\mathfrak X}\to({\mathbb C},0)$. We can suppose that the curves of the sets
${\mathcal C}^{\mathbb R}_1({\mathfrak X}_t,D,{\boldsymbol w}_t)$, $t>0$, form disjoint equisingular families.
Their limits at $t=0$ are as follows.

\begin{lemma}\label{l6}
The limit at $t=0$ of any family $C_t\in{\mathcal C}^{\mathbb R}_1({\mathfrak X}_t,D,{\boldsymbol w}_t)$, $t>0$, is
a curve $C_0=C\cup(C'\cup C'')$, where \begin{enumerate}\item[(i)] $C'\subset Y$ is a real rational curve in the linear system $|D-mE|_Y$
for some $m\ge1$, which passes through ${\boldsymbol w}'$ and transversally intersects $E$ in
$m$ distinct pairs of complex conjugate points, \item[(ii)] the curve $C'\subset Z$ is smooth rational,
representing the divisor class $-K_Z/2$, passing through $w_1$, and intersecting $E$ at some pair of complex conjugate points of $C\cap E$,
the curve $C''$ consists of $(m-1)$ pairs of complex conjugate lines
that generate the two rulings of $Z$ and pass
through $(C\cap E)\setminus(C'\cap E)$.\end{enumerate}
Furthermore, any curve $C\cup(C'\cup C'')$ as above is a limit of a unique family $C_t\in{\mathcal C}^{\mathbb R}_1({\mathfrak X}_t,D,{\boldsymbol w}_t)$, $t>0$.
\end{lemma}

{\bf Proof.} The part $C_0\cap Z$ is a nonempty real curve
passing through $w_1$. It then belongs to the linear system
$|mE|_Z$ for some $m\ge1$, and hence $C=C_0\cap Y$ belongs to $|D-mE|_Y$, $m\ge0$. Since $C\supset{\boldsymbol w}'$, the dimension count in
\cite[Proposition 2.1]{Shoval_Shustin:2013} and the genus
bound yield that either $C$ is irreducible of genus $0$ or $1$, or
$C$ consists of two components, one rational and
one elliptic. In both cases, the components of $C$ are real
and intersect $E$ in pairs of complex conjugate points.
Note that $C$ has no elliptic component. Indeed, otherwise,
the curve $C_0\cap Z$ would consists of lines from
the rulings of $Z$ and would not hit a generic point
$w_0\in{\mathbb R} Z$, since the family of real
elliptic curves in $|D-mE|_Y$ passing through ${\boldsymbol w}'$ has real dimension
one (see \cite[Proposition 2.1]{Shoval_Shustin:2013}).
Hence $C$ is real, irreducible, rational, and
intersects $E$ in $m$
distinct pairs of complex conjugate points.
The asserted structure of $C_0\cap Z$ follows immediately.

The existence and uniqueness of a family $C_t\in{\mathcal C}^{\mathbb R}_1({\mathfrak X}_t,((F_0)_t,(F_1)_t),{\boldsymbol w}_t)$, $t>0$, with
a prescribed limit $C\cup(C'\cup C'')$ satisfying conditions
(i), (ii), follow, for instance, from
\cite[Theorem 2.8]{Shustin_Tyomkin:2006}.
\proofend

Observe that the curves $C_t$ coming from a limit curve $C_0=C\cup(C'\cup C'')$ with
$C\in|D-mE|_Y$ have precisely $m-1$ solitary nodes in the component $(F_1)_t
\subset{\mathbb R}{\mathfrak X}_t$ and no other real nodes. Hence,
$$W_1(X,D,(F_0,F_1),(1,1),0)=\sum_{m\ge1}(-1)^{m-1}2^{m-1}mW(Y,D-mE,{\boldsymbol w}')\ ,$$
$$W_{1,(1,-DK_X-1)}(X,D,(F_0,F_1),(1,-1),0)=\sum_{m\ge1}2^{m-1}mW(Y,D-mE,{\boldsymbol w}')\ ,$$
where $W(Y,D-mE,{\boldsymbol w}')=\sum_C(-1)^{s(C)}$ with $C$ running over all real rational curves in the linear system
$|D-mE|_Y$ passing through ${\boldsymbol w}'$, and $s(C)$ is the total number of solitary nodes of $C$. Thus, we obtain
$$W_1(X,D,(F_0,F_1),(1,1),0)+W_{1,(1,-DK_X-1)}(X,D,(F_0,F_1),(1,-1),0)$$
$$=\sum_{m\ge1}2^{2m-1}(2m-1)W(Y,D-(2m-1)E,{\boldsymbol w}')\ .$$
On the other hand, it follows from \cite[Theorem 6(2) and Proposition 35]{Itenberg_Kharlamov_Shustin:2013} that
$$W(X,D',F_0,[F_1])=2\sum_{m\ge1}2^{2m-1}W(Y,D'-(2m-1)E,{\boldsymbol w}')$$ for any divisor class $D'\in\Pic^{\mathbb R}(X)$,
where $$W(X,D',F_0,[F_1])=\sum_{C\in{\mathcal C}^{\mathbb R}_0(X,D',{\boldsymbol w}')}(-1)^{s(C;F_0)}$$
is the (rational) Welschinger invariant (in the notation of \cite{Itenberg_Kharlamov_Shustin:2013}). So,
$$W_1(X,D,(F_0,F_1),(1,1),0)+W_{1,(-DK_X-1,1)}(X,D,(F_0,F_1),(1,-1),0)$$
\begin{equation}=\frac{1}{2}W(X,D,F_0,[F_1])+\sum_{m\ge1}W(X,D-2mE,F_0,[F_1])\ ,\label{e14}\end{equation} and we immediately
derive relations (\ref{e12}), (\ref{e13}) from the positivity and asymptotics
of Welschinger invariants $W(X,D',F_0,[F_1])$ established in
\cite[Theorem 7]{Itenberg_Kharlamov_Shustin:2013}.

\subsection{Proof of Proposition \ref{p3}}
Our argument is completely parallel to that in the proof of statement (iii) of Proposition \ref{p2}
in Section \ref{sec-iii}. First,
we construct a conjugation-invariant family ${\mathfrak X}\to({\mathbb C},0)$ of
surfaces along which the component $F_g$ (as $g=2$ or $3$) collapses, and $X$
degenerates into the
union of a real nodal del Pezzo surface and a quadric surface, intersecting along
a real rational curve $E$ with the empty real part. Then, similarly to (\ref{e14}) we derive
$$W_{2,{\underline{r}}'}(X,D,{\underline{F}}',(1,1,\pm1))=\frac{1}{2}W_1(X,D,(F_0,F_1),(1,1),0)$$ \begin{equation}+
\sum_{m\ge1}W_1(X,D-2mE,(F_0,F_1),(1,1),0)\label{e15}\end{equation} and
$$W_{3,{\underline{r}}''}(X,D,{\underline{F}}'',(1,1,\pm1,\pm1))=\frac{1}{2}W_{2,{\underline{r}}'}
(X,D,{\underline{F}}',(1,1,\pm1))$$ \begin{equation}+
\sum_{m\ge1}W_{2,{\underline{r}}'}(X,D-2mE,{\underline{F}}',(1,1,\pm1),0)\ ,\label{e16}\end{equation} provided we establish the following analog of the
vanishing statement in \cite[Proposition 35]{Itenberg_Kharlamov_Shustin:2013}:

\begin{lemma}\label{l7}
(1) Let $X,D,r_0,r_1$ be as in Proposition \ref{p3}(1). Then
$$W_1(X,D,(F_0,F_1),(1,1),[F_2])=0\ .$$

(2) Let $X,D,r_0,r_1$ be as in Proposition \ref{p3}(2). Then
$$W_2(X,D,(F_0,F_1,F_2),(1,1,\varepsilon_2),[F_3])=0,\quad\varepsilon_2=\pm1\ .$$
\end{lemma}

Observe that formula (\ref{e15}) and Proposition \ref{p2}(i,ii) yield the first statement of Proposition \ref{p3},
and subsequently formula (\ref{e16}) yields the second statement of
Proposition \ref{p3}.

{\bf Proof of Lemma \ref{l7}.} We prove the first statement; the second one can be proved in the same way.

One can check that the assumption $p_a(D)\ge2$ yields $-DK_X>2$, thus, we can assume
that $r_1>1$. As in Section \ref{sec-i}, we consider an elliptic
ABV family ${\mathfrak X}\to({\mathbb C},0)$
such that the components $F_1,F_2$ of $X={\mathfrak X}_t$ ($t>0$) degenerate into $F\cup {\mathbb R} Z$, where $Y\simeq {\mathbb R} Z\simeq S^2$,
$F\setminus{\mathbb R} E=F_+\cup F_-$, ${\mathbb R} Z\setminus{\mathbb R} E={\mathbb R} Z_+\cup{\mathbb R} Z_-$, and we suppose that the limit of $F_1$
(respectively, $F_2$) is $F_+\cup{\mathbb R} Z_+$ (respectively, $F_-\cup{\mathbb R} Z_-$). Then (for an appropriate $D_0\in\Pic^{\mathbb R}(Y)$, $D_0\ge D$) we choose a
$D_0$-CH-configuration ${\boldsymbol w}_0$ of $-DK_X$ real points
on $Y$: $r_0$ points on the component $F_0$ of ${\mathbb R} Y$ (the limit of $F_0$) and $r_1$ points in $F_+$. Similarly to Lemma \ref{l2}, we have
$$W_1(X,D,(F_0,F_1),(1,1),[F_2])=W_1(Y,D,(F_0,F_+),(1,1),{\boldsymbol w}_0)\ ,$$ where
$$W_{1,(r_0,r_1)}(Y,D,(F_0,F_+),(1,1),{\boldsymbol w}_0)=\sum_{C\in{\mathcal C}_1^{\mathbb R}(Y,D,{\boldsymbol w}_0)}(-1)^{s(C,F_0\cup F)}\ .$$
As in the proof of Lemma \ref{l3}, we specialize a suitable point $w\in {\boldsymbol w}_0\cap F_+$ to ${\mathbb R} E$, and then each curve
$C\in {\mathcal C}_1^{\mathbb R}(Y,D,{\boldsymbol w}_0)$ will degenerate into the union $C'\cup E$, where $C'\in|D-E|$ is a real immersed elliptic curve,
passing through ${\boldsymbol w}'={\boldsymbol w}_0\setminus\{w\}$,
intersecting $E$ at one point, and having there a smooth branch quadratically tangent to $E$ (the other option (b)
mentioned in the proof of Lemma \ref{l3} is not possible, since $C\cap F_-$ is finite). By \cite[Proposition 2.8(2)]{Shoval_Shustin:2013}, each curve $C'\cup E$ gives rise to two
curves in ${\mathcal C}^{\mathbb R}_1(Y,D,{\boldsymbol w}_0)$, which are distinguished by (two) deformation patterns given
in \cite[Lemma 2.10(2)]{Shoval_Shustin:2013}, and which have opposite Welschinger signs (see \cite[Proposition 6.1(i)]{Shustin:2006}), and
therefore do not contribute to
$W_{1,(r_0,r_1)}(Y,D,(F_0,F_+),(1,1),{\boldsymbol w}_0)$. \proofend

\section*{Funding}
The research has been supported by the German-Israeli Foundation, grant no. 1174-197.6/2011, and
by the Hermann-Minkowski-Minerva Center for
Geometry
at the Tel Aviv University.

{\small
\section*{Acknowledgements} The main part of this work has been performed during the author's visit to the
Centre Interfacultaire Bernoulli at the \'Ecole Polytechnique Fed\'erale de Lausanne and to High School of Economics, Moscow,
in 2014. The author is very grateful to these institutions for hospitality and excellent working conditions. Special thanks are due to the unknown referees
for numerous important remarks.}

\section{Appendix A: Degeneration and deformation of curves on rational surfaces}\label{sec-aa}

\subsection{Moduli spaces of curves}

Let
$\Sigma$ be a smooth projective rational surface
and $D\in\Pic(\Sigma)$. Denote by $\overline{\mathcal M}_{g,n}(\Sigma,D)$,
$g\ge0$,
the space of the isomorphism classes of
pairs $(\nu:\hat C\to\Sigma,{\boldsymbol p})$, where $\hat C$ is either a Riemann surface of genus $g$
or a connected reducible nodal curve of arithmetic genus $g$,
$\nu_*\hat C\in|D|$,
${\boldsymbol p}=(p_1,...,p_n)$ is a sequence of distinct smooth points of $\hat C$, and each component
$C'$ of $\hat C$ of genus $g'$, which is contracted
by $\nu$,
contains at
least $3-2g'$ special points. This
moduli space is a projective scheme (see \cite{Fulton_Pandharipande:1997}),
and
there are natural morphisms
$$\Phi_{\Sigma,D}:\overline{\mathcal M}_{g,n}(\Sigma,D)\to|D|,\quad [\nu:\hat
C\to\Sigma,{\boldsymbol p}]\mapsto\nu_*\hat C\ ,$$ $$\Ev:\overline{\mathcal
M}_{g,n}(\Sigma,D)\to\Sigma^n,\quad [\nu:\hat C\to\Sigma,{\boldsymbol p}]\mapsto\nu({\boldsymbol p})\
.$$ For any subscheme ${\mathcal
V}\subset\overline{\mathcal M}_{0,n}(\Sigma,
D)$, define the {\it intersection dimension} $\idim\mathcal V$ of $\mathcal V$ as follows:
$$\idim{\mathcal V}=\dim(\Phi_{\Sigma,D}\times\Ev)({\mathcal V})\ ,$$ where the latter value is the maximum over the dimensions of
all irreducible components.

Put
\begin{eqnarray}
{\mathcal M}_{g,n}^{br}(\Sigma,D)&=&\{[\nu:\hat C\to\Sigma,{\boldsymbol p}]\ \in {\mathcal M}_{g,n}(\Sigma,D)\ |\ \hat C\ \text{is smooth, and}\nonumber\\
& &\qquad\qquad\qquad\qquad\qquad \nu\ \text{is birational on to}\ \nu(\hat C)\},\nonumber\\ {\mathcal M}_{g,n}^{im}(\Sigma,D)&=&\{[\nu:\hat C\to\Sigma,{\boldsymbol p}]
\in{\mathcal M}_{g,n}^{br}(\Sigma,D)\ |\ \nu\ \text{is an immersion}\}\ .\nonumber
\end{eqnarray} Denote by $\overline{{\mathcal M}^{br}_{g,n}}(\Sigma,D)$ the
closure of ${\mathcal M}^{br}_{g,n}(\Sigma,D)$ in
$\overline{\mathcal M}_{g,n}(\Sigma,D)$, and introduce also the space
$${\mathcal M}'_{g,n}(\Sigma,D)=\{[\nu:\hat C\to\Sigma,{\boldsymbol p}]
\in \overline{{\mathcal M}^{br}_{g,n}}(\Sigma,D)\ |\ \hat C\ \text{is smooth}\}\ .$$

\begin{lemma}\label{ln15} For any element $[\nu]=[\nu:\hat
C\to\Sigma,{\boldsymbol p}]\in{\mathcal M}_{g,n}^{br}(\Sigma,D)$, the map $\Phi_{\Sigma,D}\times\Ev$ is
injective
in a neighborhood $[\nu]$, and, for the germ at $[\nu]$ of any irreducible subscheme ${\mathcal V}\subset
{\mathcal M}_{0,n}^{br}(\Sigma,D)$, we have
$$\dim {\mathcal V}=\idim {\mathcal V}\ .$$
\end{lemma}

\subsection{Curves on
del Pezzo and uninodal del Pezzo surfaces}\label{sec-dpudp}

Let $\Sigma$ be the plane ${\mathbb P}^2$ blown up at eight distinct points. Denote
by ${\mathcal D}$ the Kodaira-Spencer-Kuranishi space of all complex structures on
the smooth four-fold $\Sigma$, factorized by the action of diffeomorphisms homotopic to the identity.
It contains an open dense subset ${\mathcal D}^\DP$ consisting of del Pezzo surfaces (of degree $1$),
that is, surfaces with an ample effective anticanonical class. We call a surface $Y\in{\mathcal D}$
uninodal del Pezzo,
if it contains a smooth rational $(-2)$-curve $E_Y$, and $-K_Y C>0$
for each irreducible curve $C\ne E_Y$ (in particular, $C^2\ge-1$). Denote by
${\mathcal D}^\DP(A_1)\subset{\mathcal D}$ the subspace formed by uninodal del Pezzo surfaces.
Observe that ${\mathcal D}^\DP(A_1)$ has codimension $1$ in ${\mathcal D}$, and
${\mathcal D}\setminus({\mathcal D}^\DP\cup{\mathcal D}^\DP(A_1))$ is of codimension $\ge2$ in
${\mathcal D}$.

Through all this section we use the notation
$$n=-DK_\Sigma+g-1\ .$$

\begin{lemma}\label{ln11}
If $\Sigma$ is a smooth rational surface and
$-DK_{\Sigma}>0$, then
the space ${\mathcal M}^{im}_{g,0}(\Sigma,D)$ is either empty, or is a smooth variety
of dimension $n$.
\end{lemma}

{\bf Proof.}
Let
$[\nu:\hat C\to\Sigma]\in{\mathcal M}_{g,0}^{im}(\Sigma,D)$. The Zariski
tangent space to ${\mathcal M}_{g,0}^{im}(\Sigma,D)$ at
$[\nu]$
can be identified with
$H^0(\hat C,{\mathcal N}^{\nu}_{\hat C})$. Since
\begin{equation}\deg{\mathcal
N}^{\nu}_{\hat C}=-DK_\Sigma+2g-2>2g-2\
,\label{en7}\end{equation} we have
\begin{equation}h^1(\hat C,{\mathcal N}^{\nu}_{\hat C})=0\ ,\label{e20}\end{equation}
and hence
${\mathcal M}_{g,0}^{im}(\Sigma,D)$ is smooth at $[\nu]$
and is
of
dimension
\begin{equation}h^0(\hat C,{\mathcal N}^{\nu}_{\hat C})=\deg{\mathcal
N}^{\nu}_{\hat C}-g+1=-DK_\Sigma+g-1=n\ .\label{en8}\end{equation}
\proofend

\begin{lemma}\label{ln12}
(1) Let $\Sigma\in{\mathcal D}^{\DP}$ and $-DK_{\Sigma}>0$. Then, the following holds:
\begin{enumerate}\item[(i)]
The space ${\mathcal M}_{g,0}^{br}(\Sigma,D)$ is either empty or satisfies
$\dim{\mathcal M}_{g,0}^{br}(\Sigma,D)\le n$.
\item[(ii)] If
either $g > 0$ or $D\ne-K_\Sigma$, then ${\mathcal
M}_{g,0}^{im}(\Sigma,D)\subset{\mathcal M}_{g,0}^{br,n}(\Sigma,D)$ is
an
open dense subset,
where ${\mathcal M}_{g,0}^{br,n}(\Sigma,D)$ denotes the union of the components of
${\mathcal M}_{g,0}^{br}(\Sigma,D)$ of dimension $n$. \item[(iii)]
There
exists
an open dense
subset $U^\DP
\subset{\mathcal D}^{\DP}$ such that, if
$\Sigma \in U^\DP$,
then ${\mathcal M}_{0,0}(\Sigma,-K_\Sigma)$ consists of
twelve elements, each corresponding to a rational nodal curve.\end{enumerate}

(2) There exists an open dense subset $U^\DP(A_1)\subset{\mathcal D}(A_1)$ such that
if $\Sigma\in U^\DP(A_1)$ and $-DK_{\Sigma}>0$,
then
\begin{enumerate}
\item[(i)] $\idim{\mathcal M}'_{g,0}(\Sigma,D)\le n$;
\item[(ii)] a generic
element $[\nu:\hat C\to\Sigma]$ of any irreducible component ${\mathcal V}$ of
${\mathcal M}'_{g,0}(\Sigma,D)$, such that $\idim{\mathcal V}=n$, is an immersion, and
the divisor
$\nu^*(E_\Sigma)$
consists of $DE_\Sigma$ distinct points.
\end{enumerate}
\end{lemma}

{\bf Proof.}
Let $\Sigma\in{\mathcal D}^{\DP}\cup{\mathcal D}^\DP(A_1)$. All the statements for the case of an effective
$-K_\Sigma-D$ immediately follow from elementary properties of plane lines, conics, and cubics.
In particular, a general element of ${\mathcal D}^\DP\setminus U^\DP$ is the plane
blown
up
at eight generic points
on a cuspidal cubic.
So, in the sequel we suppose that $-K_\Sigma-D$ is not effective.

In view of Lemma \ref{ln11}, to complete the proof of
statements (1) and (2i) it is enough to show that $$\dim ({\mathcal M}_{g,0}^{br,n}(\Sigma,D)\setminus{\mathcal
M}_{g,0}^{im}(\Sigma,D))<n\quad\text{and}\quad\dim({\mathcal M}'_{g,0}
(\Sigma,D)\setminus{\mathcal M}^{br,n}_{g,0}(\Sigma,D))<n\ .$$

Note, first, that, in the case $n=0$, we have $g=0$ and $-DK_\Sigma=1$, and the
curves $C\in\Phi_{\Sigma,D}({\mathcal M}_{g,0}^{br,0}(\Sigma,D))$ are
nonsingular due to the bound
\begin{equation}-DK_\Sigma\ge (C\cdot C')(z)\ge s\ ,\label{en6}\end{equation} coming from the intersection of $C$ with
a curve $C'\in|-K_\Sigma|$ passing through a point $z\in C$, where $C$
has multiplicity $s$. Thus, further on we suppose that $n>0$.

Let ${\mathcal V}_2$ be an irreducible component
of ${\mathcal M}_{g,0}^{br,n}(\Sigma,D)\setminus {\mathcal
M}_0^{im}(\Sigma,D)$, $[\nu:\hat C\to\Sigma]\in{\mathcal V}_2$ a generic
element, and let $\nu$ have $s\ge1$ critical points
of multiplicities $m_1\ge...\ge m_s\ge2$. Particularly, bound (\ref{en6}) gives
\begin{equation}-DK_\Sigma\ge m_1\ .\label{en17}\end{equation}
Then ({\it cf.} \cite[First formula in the proof of Corollary
2.4]{Caporaso_Harris:1998}),
$$\dim{\mathcal
V}_2\le h^0(\hat C,{\mathcal
N}^{\nu}_{\hat C}/\Tors({\mathcal N}^{\nu}_{\hat C}))\ ,$$ where the
normal sheaf ${\mathcal N}^{\nu}_{\hat C}$ on $\hat C$ is defined as
the cokernel of the map $d\nu:{\mathcal T}\hat C\to\nu^*{\mathcal
T}\Sigma$, and $\Tors(*)$ is the torsion sheaf.
It follows from \cite[Lemma 2.6]{Caporaso_Harris:1998} ({\it cf.} also the computation in \cite[Page 363]{Caporaso_Harris:1998})
that $\deg\Tors({\mathcal
N}^{\nu}_{\hat C})=\sum_i(m_i-1)$, and hence
\begin{equation}\deg{\mathcal N}^{\nu}_{\hat C}/\Tors({\mathcal
N}^{\nu}_{\hat C})=-DK_\Sigma+2g-2-\sum_{i=1}^s(m_i-1)\label{en24}\end{equation}
which yields
\begin{eqnarray}\dim{\mathcal V}_2
&\le& h^0(\hat C,{\mathcal
N}^{\nu}_{\hat C}/\Tors({\mathcal
N}^{\nu}_{\hat C}))\nonumber\\
&=&\max\{\deg{\mathcal N}^{\nu}_{\hat C}/\Tors({\mathcal
N}^{\nu}_{\hat C})-g+1,\
g\}\overset{\text{(\ref{en17})}}{\le}n-(m_1-1)<n,\label{en5}\end{eqnarray}

Let ${\mathcal V}$ be an irreducible component
of ${\mathcal M}'_{g,0}(\Sigma,D)\setminus{\mathcal M}_{g,0}^{br,n}(\Sigma,D)$.
Then a generic element $[\nu:\hat C\to\Sigma]\in{\mathcal V}$ satisfies $\nu_*\hat C=sC$ for some $s\ge2$ and
some reduced, irreducible curve $C\subset\Sigma$. It follows from
the Riemann-Hirwitz formula, that
$g(C)-1\le\frac{1}{s}(g-1)$, and hence
$$\idim{\mathcal
V}\le-CK_\Sigma+g(C)-1\le-\frac{1}{s}(DK_\Sigma+g-1)<-DK_\Sigma+g-1=n\ .$$

To complete the proof of (2ii), let us assume that $\dim{\mathcal V}=r$ and the divisor $\nu^*(E_\Sigma)$
contains a multiple point $sz$, $s\ge2$.
In view of $DE_\Sigma \ge s$ and $(-K_\Sigma - E_\Sigma)D\ge0$
(remind that $D$ is irreducible and $-K_\Sigma-D$ is not effective), we have
$-DK_\Sigma\ge s$. Furthermore,  $T_{[\nu]}{\mathcal V}$ can be identified with
a subspace of $H^0(\hat C,{\mathcal N}_{\hat C}^{\nu}(-(s-1)z))$ ({\it cf.} \cite[Remark in page 364]{Caporaso_Harris:1998}). Since
$$\deg {\mathcal N}_{\hat C}^{\nu}(-(s-1)z))=-DK_\Sigma+2g-1-s
\overset{-DK_\Sigma\ge s}{>}2g-2\ ,$$
we have
$$H^1(\hat C,{\mathcal N}_{\hat C}^{\nu}(-(s-1)z^*))=0\ ,$$ and hence $$
\dim{\mathcal V}\le h^0(\hat C,{\mathcal N}_{\hat C}^{\nu}(-(s-1)z^*))=n-(s-1)<n$$ contrary to the assumption
$\dim{\mathcal V}=n$. \proofend

\begin{lemma}\label{ln17}
There exists an open dense subset $V^{\DP}\subset{\mathcal D}^{\DP}$
such that, for each $\Sigma\in V^{\DP}$, the set
of effective divisor classes $D\in\Pic(\Sigma)$ satisfying $-DK_\Sigma=1$ is finite, the set of rational curves
in the corresponding linear systems $|D|$ is finite, and any two such rational curves $C_1,C_2$ either coincide, or are disjoint, or
intersect in $C_1C_2$ distinct points.
\end{lemma}

{\bf Proof}.
For the proof see \cite[Lemma 10]{Itenberg_Kharlamov_Shustin:2012}.
\proofend

\begin{lemma}\label{ln16}
For each surface $\Sigma\in U^\DP\cap V^{\DP}$ ,
each divisor class $D \in \Pic(\Sigma)$ with
$-DK_\Sigma>0$ and $D^2\ge -1$,
and each irreducible component ${\mathcal V}$ of
$\overline{{\mathcal M}^{br,n}_{g,0}}(\Sigma,D)\setminus{\mathcal M}_{g,0}^{br,n}(\Sigma,D)$
with $\idim{\mathcal V}=n-1$,
a generic element $[\nu:\hat C\to\Sigma]\in{\mathcal V}$ is such that
\begin{enumerate}\item[(i)] either
$\hat C=\hat C_1\cup\hat C_2$ with $\hat C_1,\hat C_2$ smooth Riemann surfaces of genera $g_1,g_2$, respectively, such that $g=g_1+g_2$;
furthermore, $|\hat C_1\cap\hat C_2|=1$,
$[\nu|_{\hat C_i}:\hat C_i\to\Sigma]\in{\mathcal M}_{g_i,0}^{im}(\Sigma,D_i)$,
where $C_1=\nu(\hat C_1)\ne C_2=\nu(\hat C_2)$, $D_1D_2>0$, and $-D_iK_\Sigma>0$, $D_i^2\ge -1$
for each $i=1,2$, and, in addition, at any point $z\in C_1\cap C_2$, any
component of $(C_1,z)$ intersects any component of $(C_2,z)$ transversally;
\item[(ii)] or $D=-2K_\Sigma$, $g=0$, $\hat C=\hat C_1\cup\hat C_2$,
$|\hat C_1\cap\hat C_2|=1$,
$\nu|_{\hat C_1}$ and $\nu|_{\hat C_2}$ are immersions of $\hat C_1\simeq\hat C_2\simeq{\mathbb P}^1$
on to the same uninodal curve $C\in|-K_\Sigma|$;
\item[(iii)] or $D=-2K_\Sigma$, $\hat C$ is a smooth elliptic curve, $\nu:\hat C\to C=\nu(\hat C)$
is an unramified double covering.
\end{enumerate}
Furthermore,
$\nu$ is always an immersion ({\it i.e.}, a local isomorphism on to the image), and the germ of $\overline{{\mathcal
M}^{br,n}_{g,0}}(\Sigma,D)$ at $[\nu]$ is smooth.
\end{lemma}

{\bf Proof.} Let ${\mathcal V}$ be an irreducible component of
$({\mathcal M}'_{g,0}(\Sigma,D)\cap\overline{{\mathcal M}^{br,n}_{g,0}}
(\Sigma,D))\setminus{\mathcal M}^{br,n}_{g,0}(\Sigma,D)$ such that $\idim{\mathcal V}=n-1$
($\idim{\mathcal V}$ cannot be bigger by Lemma
\ref{ln12}(i)).
Then its generic element $[\nu:\hat C\to\Sigma]$ is such that $\nu_*\hat C=sC$ with a reduced, irreducible $C$, $s\ge2$. By the Riemann-Hurwitz formula,
$g-1=s(g(C)-1)+\rho/2$, where $\rho$ is the total ramification index of the map
$\nu^\vee:\hat C\to C^\vee$, $C^\vee$ being the normalization of $C$.
By Lemma \ref{ln12}(i), $$\idim{\mathcal V}=n-1=-sCK_\Sigma+g-2\le-CK_\Sigma+g(C)-1\ ,$$
which together with the above Riemann-Hurwitz formula yields
$$(s-1)(-CK_\Sigma+g(C)-1)+\frac{\rho}{2}\le1\ .$$ It follows that
\begin{itemize}\item either $s=2$, $g=g(C)=1$,
$-CK_\Sigma=1$, $\rho=0$, and hence $C\in|-K_\Sigma|$, which meets one of the cases in
statement (i);
\item or $s=2$, $-CK_\Sigma=1$, and $g(C)=0$, which yields $g=0$ and, in view of the adjunction formula,
$C^2=-1$, or $C^2\ge1$; both cases are not possible: the former one is excluded by the assumption $D^2\ge-1$,
whereas the latter one leaves the only option  of $C\in|-K_\Sigma|$ a
uninodal curve, however, in such a case
the map $\nu$ cannot be deformed into an element of ${\mathcal M}^{br}_{0,0}(\Sigma,-2K_\Sigma)$, since
the
deformed map would birationally send ${\mathbb P}^1$ on to a curve with
$\delta$-invariant $\ge 4$ in a neighborhood of the node of $C$,
which is bigger
than the arithmetic genus, $((-2K_\Sigma)^2+(-2K_\Sigma)K_\Sigma)/2+1=2$.
\end{itemize}

Now let $[\nu:\hat C\to\Sigma]$ be a generic element of an irreducible component
${\mathcal V}$ of $\overline{{\mathcal M}^{br,n}_{g,0}}(\Sigma,D)\setminus{\mathcal
M}'_{g,0}(\Sigma,D)$ with $\idim{\mathcal V}=n-1$. Then $\hat C$ has $s\ge2$ components
$\hat C_1,...,\hat C_s$ of genera $g_1,...,g_s$, respectively, and $l\ge s-1$ nodes.
It follows that $g=g_1+...+g_s+l-s+1$, and, by Lemma \ref{ln12}(i),
$$\idim{\mathcal V}=n-1=-DK_\Sigma+g-2\le-DK_\Sigma+(g_1+...+g_s)-s\ ,$$ and hence $l=1$,
$s=2$, $g=g_1+g_2$. By Lemma \ref{ln12}(ii), both $\nu|_{\hat C_1}$ and $\nu|_{\hat C_2}$ are immersions.
Note that the case $\nu(\hat C_1)=\nu(\hat C_2)$ is possible only when $D_1=D_2$, $g_1=g_2$, and
$-D_1K_\Sigma+g_1-1=0$. Since $D_1^2=D_2^2\ge1$ in view of the adjunction formula and the
condition $D^2\ge-1$, we are left with the case $D_1=D_2=-K_\Sigma$, and $C=\nu(\hat C_1)=\nu(\hat C_2)\in|-K_\Sigma|$
a rational curve with the unique node $z$.
The map $\nu$ takes the germ $(\hat C,\hat z)$ isomorphically on to the germ
$(C,z)$, since, otherwise, we would get a deformed map $\nu$ with
the image whose
$\delta$-invariant $\ge 4$,
which is bigger
than its arithmetic genus, $((-2K_\Sigma)^2+(-2K_\Sigma)K_\Sigma)/2+1=2$.
Suppose now that $C_1=\nu(\hat C_1)\ne C_2=\nu(\hat C_2)$. Let us show that
$C_1$ and $C_2$ intersect transversally as proclaimed in statement (i), which would imply that
$\nu$ is in immersion.
If $-D_1K_\Sigma+g_1-1=-D_2K_\Sigma+g_2-1=0$, then $g_1=g_2=0$ and $-D_1K_\Sigma=-D_2K_\Sigma=1$,
which allows one to apply Lemma \ref{ln17}. If $-D_1K_\Sigma+g_1-1>0$, then we can vary
$\nu|_{\hat C_1}$ in ${\mathcal M}_{g_1,0}^{im}(\Sigma,D_1)$ and achieve the required transversality
as we did in the proof of Lemma \ref{ln12}(2ii).

At last, the proof of the smoothness of the germ of $\overline{{\mathcal
M}^{br,n}_{g,0}}(\Sigma,D)$ at $[\nu]$ literally coincides with that in \cite[Lemma 11]{Itenberg_Kharlamov_Shustin:2012}.
\proofend

\begin{lemma}\label{p2X} Let $\Sigma\in U^\DP$, $g\ge0$, and $D\in\Pic(\Sigma)$ be
an effective divisor class such that
$n=-DK_\Sigma+g-1\ge1$.
Let
${\boldsymbol w}=(w_1,...,w_n)$ be a sequence of $n$ distinct points in $\Sigma$, let
$\sigma_i$
be
smooth curve germs in $\Sigma$ centered at $w_i$, $n'<i\le
n$, for some $n'<n$, ${\boldsymbol w}'=(w_i)_{1\le i\le n'}$, and let
\begin{eqnarray}\overline{{\mathcal M}_{g,n}^{br}}(\Sigma,D;{\boldsymbol w}',\{\sigma_i\}_{n'<i\le
n})&=&\{[\nu:\hat C\to\Sigma,{\boldsymbol p}]\in\overline{{\mathcal
M}_{g,n}^{br}}(\Sigma,D)\ : \nonumber\\
& & \quad \nu(p_i)=w_i \ \text{for}\, \ 1\le i\le n',\ \nu(p_i)\in \sigma_i, \,\text{for}\, \ n'<i\le
n\}\ .\nonumber\end{eqnarray}

(1) Suppose that $[\nu:\hat C\to\Sigma,{\boldsymbol p}]$ either belongs to $\overline{{\mathcal
M}_{g,n}^{br}}(\Sigma,D;{\boldsymbol w})\cap{\mathcal
M}_{g,n}^{im}(\Sigma,D)$, or is as in Lemma \ref{ln12}(iii). If
\begin{equation}H^1(\hat C,{\mathcal N}^\nu_{\hat C}(-{\boldsymbol p}))=0\ ,
\label{eelra1}\end{equation} then $\Ev$ sends the germ of $\overline{{\mathcal
M}_{g,n}^{br}}(\Sigma,D;{\boldsymbol w}',\{\sigma_i\}_{n'<i\le n})$ at
$[\nu:{\mathbb P}^1\to\Sigma,{\boldsymbol p}]$ diffeomorphically on to $\prod_{n'<i\le
n}\sigma_i$.

(2) Suppose that $[\nu:\hat C\to\Sigma,{\boldsymbol p}]\in\overline{{\mathcal
M}_{g,n}^{br}}(\Sigma,D;{\boldsymbol w})$ is such that
\begin{itemize}\item $[\nu:\hat C\to\Sigma]\in\overline{{\mathcal
M}_{g,0}^{br}}(\Sigma,D)$ is as in Lemma \ref{ln16}(i),
\item $n'\ge-D_1K_\Sigma+g_1-1$, $\#({\boldsymbol p}\cap \hat
C_1)=-D_1K_\Sigma+g_1-1$, $\#({\boldsymbol p}\cap \hat C_2)=-D_2K_\Sigma+g_2$, the point sequences
$(w_i)_{1\le i<-D_1K_\Sigma}$ and $(w_i)_{-D_1K_\Sigma\le i\le n}$ are
generic on the curves $C_1=\nu_*\hat C_1$ and $C_2=\nu_*\hat
C_2$, respectively, and the germs $\sigma_i$, $n'<i\le n$, cross $C_2$
transversally.\end{itemize} Then $\Ev$ sends the germ of
$\overline{{\mathcal
M}_{g,n}^{br}}(\Sigma,D;{\boldsymbol w}',\{\sigma_i\}_{n'<i\le n})$ at
$[\nu:\hat C\to\Sigma,{\boldsymbol p}]$ diffeomorphically on to $\prod_{n'<i\le
n}\sigma_i$.
\end{lemma}

{\bf Proof.} The first statement follows from the fact that $\Ev$ diffeomorphically
sends the (smooth by Lemma \ref{ln11})
germ of $\overline{{\mathcal M}_{g,n}^{br}}(\Sigma,D)$ at $[\nu:\hat C\to\Sigma,
{\boldsymbol p}]$ on to the germ
of $\Sigma^r$ at ${\boldsymbol w}$.
In view of $\dim\overline{{\mathcal M}_{g,n}^{br}}(\Sigma,D)=2n$
(see
Lemma \ref{ln12}(i)) it is
sufficient
to show that
the Zariski tangent space to $\Ev^{-1}({\boldsymbol w})$ is zero-dimensional,
which is
equivalent to
\begin{equation}h^0(\hat C,{\mathcal N}_{\hat C}^\nu(-{\boldsymbol p}))=0
\label{eK5}\end{equation} that in turn immediately follows from (\ref{en8}) and
(\ref{eelra1}).

In the second case, by Lemma \ref{ln16}, the germ of $\overline{{\mathcal
M}_{g,n}^{br}}(\Sigma,D)$ at
$[\nu:\hat C\to\Sigma,{\boldsymbol p}]$ is smooth. The general position of the points ${\boldsymbol w}$ on the
curve $C_1\cup C_2$ yields (\ref{eelra1}), which similarly to the preceding paragraph
suffices for the required diffeomorphism ({\it cf.} proof of \cite[Lemma
12(2)]{Itenberg_Kharlamov_Shustin:2012}). \proofend

Consider a proper submersion $\widetilde\Sigma\to({\mathbb C},0)$ a smooth three-fold $\widetilde\Sigma$ such that
$\Sigma=\widetilde\Sigma_0\in U^\DP(A_1)$ and $\widetilde\Sigma_t\in U^\DP$ for all $t\ne 0$. Choose a divisor class $D\in\Pic(\Sigma)$
such that $-DK_\Sigma>0$ and a nonnegative integer $g$. Let ${\boldsymbol w}_t\in\widetilde\Sigma_t^n$, $t\in({\mathbb C},0)$, be a smooth family of
configurations of distinct points such that ${\boldsymbol w}={\boldsymbol w}_0$ is disjoint with the
$(-2)$-curve $E_\Sigma\subset\Sigma$.

\begin{lemma}\label{lem-abv}
There exists an open dense subset ${\mathcal U}_n\subset\Sigma^n$ such that, if ${\boldsymbol w}\in
{\mathcal U}_n$, then
\begin{enumerate}\item[(i)] for any $m\ge0$ and any element $[\nu:\hat C\to\Sigma,{\boldsymbol p}]\in
{\mathcal M}'_{g,n}(\Sigma,D-mE_\Sigma)$ with $\nu({\boldsymbol p})={\boldsymbol w}$, the map $\nu$ is an immersion,
and the divisor $\nu^*(E_\Sigma)\subset\hat C$ consists of $DE_\Sigma+2m$ distinct points;
\item[(ii)] each element $[\nu:\hat C\cup(\hat E_1\cup...\cup\hat E_m)\to\Sigma,{\boldsymbol p}]
\in\overline{\mathcal M}_{g,n}(\Sigma,D)$ such that
\begin{itemize}\item $[\nu:\hat C\to\Sigma,{\boldsymbol p}]\in{\mathcal M}'_{g,n}(\Sigma,D-mE_\Sigma)$, $\nu({\boldsymbol p})
={\boldsymbol w}$,
\item $\hat E_1\simeq...\simeq\hat E_m\simeq{\mathbb P}^1$ and $\nu$ takes each of $\hat E_1,...,\hat E_m$ isomorphically on
to $E_\Sigma$,
\item $\hat E_i\cap\hat E_j=\emptyset$ as $i\ne j$, and $|\hat C\cap\hat E_i|=1$ for all $i=1,...,m$,
\end{itemize} admits an extension to a smooth family $[\nu_t:\hat C_t\to\widetilde\Sigma_t,{\boldsymbol p}_t]
\in\overline{\mathcal M}_{g,n}(\widetilde\Sigma_t,D)$, $t\in({\mathbb C},0)$,
where $\nu_t({\boldsymbol p}_t)={\boldsymbol w}_t$ and $[\nu_t:\hat C_t\to\widetilde\Sigma_t,{\boldsymbol p}_t]
\in{\mathcal M}^{im}_{g,n}(\widetilde\Sigma_t,D)$,
\item[(iii)]the set of families introduced in item (ii) is in one-to-one correspondence with
each of the sets ${\mathcal C}_g(\widetilde\Sigma_t,D,{\boldsymbol w}_t)$, $t\ne0$.
\end{enumerate}
\end{lemma}

{\bf Proof.}
The statement follows from Lemma \ref{ln12}(2) and \cite[Theorem 4.2]{Vakil:2000}.
\proofend

\subsection{Deformation of isolated curve singularities}

\subsubsection{Local equigeneric deformations}
Let $\Sigma$ be a smooth algebraic surface, and
$z$ be
an isolated singular point of a curve $C\subset\Sigma$.
Denote by $J_{C,z}\subset{\mathcal O}_{C,z}$ the Jacobian ideal,
and by $J^{cond}_{C,z}\subset{\mathcal O}_{C,z}$ the local conductor ideal, defined as
$\Ann{\mathcal O}_{C^\vee}/{\mathcal O}_{C,z}$, where $C^\vee\to(C,z)$ is the normalization.
If $f(x,y)=0$ is an equation of $(C,z)$ in some local coordinates $x,y$ in $(\Sigma,z)$, then
$$J_{C,z}=\langle f_x,f_y\rangle,\quad
J^{cond}_{C,z}=\{g\in{\mathcal O}_{C,z}\ :\
\ord(g|_{C_i})\ge\ord(f'|_{C_i})-\mt(C_i)+1,\ i=1,...,m\}\ ,$$
where $C_1,...,c_m$ are all the components of $(C,z)$, $f'=\alpha f_x+\beta f_y$ a generic polar,
and $\mt(C_i)$ is
the intersection number of $C_i$ with a generic smooth line through $z$ ({\it cf.} \cite[Section 4.2.4]{Dolgachev:2013}).

Let $B_{C,z}$ be
the base of a semiuniversal
deformation
of the germ $(C,z)$.
This base
can
be
identified with ${\mathcal O}_{C,z}/J_{C,z}\simeq{\mathbb C}^{\tau(C,z)}$, where
$J_{C,z}\subset{\mathcal O}_{C,z}$
is the Jacobian ideal, $\tau(C,z)$ the Tjurina number.

Denote by
$B^{\,eg}_{C,z} \subset B_{C,z}$ the equigeneric locus that
parameterizes local
deformations of $(C,z)$ with the constant
$\delta$-invariant equal to $\delta(C,z)$. The following lemma presents the properties of
$B^{\,eg}_{C,z}$, which we will need.

\begin{lemma}\label{eg1}
The locus $B^{\,eg}_{C,z}$ is irreducible
and has codimension
$\delta(C,z)$ in $B_{C,
z}$.
The subset $B^{\,eg,im}_{C,z} \subset B^{\,eg}_{C,z}$
that parameterizes the immersed deformations is open and dense in  $B^{\,eg}_{C
, z}$,
and consists only of smooth points of $B^{\,eg}_{C, z}$.
The subset $B^{\,eg,nod}_{C,z} \subset B^{\,eg}_{C,z}$ that
parameterizes the
nodal deformations is also open and dense.
The complement $B^{\,eg}_{C,z}\setminus B^{\,eg,nod}_{C,z}$ is the closure of
three codimension-one strata: $B^{\,eg}_{C,z}(A_2)$
that parameterizes deformations with one cusp
$A_2$ and $\delta(C,z)-1$ nodes,
$B^{\,eg}_{C,z}(A_3)$
that parameterizes
deformations with one tacnode $A_3$ and $\delta(C,z)-2$ nodes,
and $B^{\,eg}_{C,
z}(D_4)$
that parameterizes deformations with one ordinary triple point
$D_4$ and $\delta(C,z)-3$ nodes.
The tangent cone
$T_0B^{\,eg}_{C,z}$ (defined as the limit of the
tangent spaces at points of $B^{\,eg,im}_{C,z}$) can be
identified with $J^{cond}_{C,z}/J_{C,z}$.
\end{lemma}

{\bf Proof.} The statement follows from \cite[Item (iii) in page 435, Theorem 1.4,
Theorem 4.15, and Proposition 4.17]{Diaz_Harris:1984}.\proofend

\subsubsection{Local invariance of Welschinger numbers}

Suppose now that $\Sigma$ possesses a real structure, $C$ is
a real curve, and $z$ is its real singular point.
Let $b\in B^{\,eg,im}_{C,z}$ be a real point, and let $C_b$
be the corresponding fiber
of the semiuniversal deformation of the germ $(C,z)$.
Choose a real point $b'\in B^{\,eg,nod}_{C,z}$ sufficiently close to $b$ and
define
Welschinger signs $$W^+_b=(-1)^{s_+(C_{b'})},\quad W^-_b=(-1)^{s_-(C_{b'})}\ ,$$ where
$s_+(C_{b'})$ (respectively, $s_-(C_{b'})$) is the number of solitary (respectively, non-solitary)
nodes of $C_{b'}$.

\begin{lemma}\label{new_lemma1}
Welschinger signs $W^+_b$ and $W^-_b$ do not depend on the choice of a real point
$b'\in B^{\,eg,nod}_{C,z}$ sufficiently close to $b$.
\end{lemma}

{\bf Proof.} Straightforward.
\proofend

\begin{lemma}\label{ln2}
Let $L_t$, $t\in(-\varepsilon,\varepsilon)\subset{\mathbb R}$, be a
continuous one-parameter family of conjugation-invariant affine
subspaces of $B_{C,z}$ of dimension $\delta(C,z)$ such that
\begin{itemize}\item $L_0$ passes through the origin and is
transversal to $T_0B^{\,eg}_{C,z}$, \item $L_t\cap B^{\,eg}_{C,z}\subset
B^{\,eg,im}_{C,z}$ for each $t \in (-\varepsilon, \varepsilon)
\setminus \{0\}$.
\end{itemize} Then,
\begin{enumerate}
\item[(i)]
the intersection $L_t\cap B^{\,eg}_{C,z}$ is finite for each $t
\in (-\varepsilon', \varepsilon') \setminus \{0\}$,
where $\varepsilon' > 0$ is sufficiently small.
\item[(ii)] the functions $W^\pm(t)=\sum_{b\in L_t\cap{\mathbb R}
B^{\,eg}_{C,z}}W^\pm_b$ are
constant in $(-\varepsilon', \varepsilon') \setminus \{0\}$,
where $\varepsilon' > 0$ is sufficiently small.
\end{enumerate}
\end{lemma}

{\bf Proof.} The statement follows from \cite[Lemma 15]{Itenberg_Kharlamov_Shustin:2012}.
\proofend

\subsubsection{Global transversality conditions}

If $C \subset \Sigma$ is a curve with isolated singularities,
we consider the joint semiuniversal deformation
for all singular points of $C$. The
base of this deformation,
the equigeneric locus,
and the tangent cone to this locus at the point corresponding
to $C$ are as follows:
$$
B_C=\prod_{z\in\Sing(C)}B_{C,z},\quad
B^{\,eg}_C=\prod_{z \in \Sing(C)}B^{\,eg}_{C,z}, \quad
T_0B^{\,eg}_C=\prod_{z \in \Sing(C)}T_0B^{\,eg}_{C,z}\ .$$

\begin{lemma}\label{leg} Let $[\nu:\hat C_1\to\Sigma,{\boldsymbol p}]\in
{\mathcal M}^{br}_{g,0}(\Sigma,D)$ and $C=\nu(\hat C)$.
Assume that
$n=-DK_\Sigma+g-1>0$. There exists an open dense subset ${\mathcal U}\subset C^n$ such that
any ${\boldsymbol p}\in{\mathcal U}$ consists of $n$ distinct points, the
image
${\boldsymbol w}=\nu({\boldsymbol p})$ is an $n$-tuple of distinct
nonsingular points of $C$, and \begin{equation}
H^0(C,{\mathcal J}^{cond}_C(-{\boldsymbol w})\otimes
{\mathcal O}_\Sigma(D))=0\ .\label{eeg}\end{equation}
Let ${\boldsymbol w}$ satisfy (\ref{eeg}), $|D|_{\boldsymbol w}\subset|D|$
be the linear subsystem of curves passing through
${\boldsymbol w}$, and $\Lambda({\boldsymbol w})\subset B_C$
be the natural image of $|D|_{\boldsymbol w}$.
\begin{enumerate}
\item[(1)] One has $\codim_{B_C}\Lambda({\boldsymbol w})=\dim
B^{\,eg}_C$, and $\Lambda({\boldsymbol w})$ intersects $T_0B^{\,eg}_C$
transversally.
\item[(2)] For any $n$-tuple ${\boldsymbol w}' \in
\Sigma^n$ sufficiently close to ${\boldsymbol w}$ and such that
$\Lambda({\boldsymbol w}')$ intersects $B^{\,eg}_C$
transversally and only at smooth points,
the natural map from the germ of ${\mathcal M}_{g,r}(\Sigma,D)$ at $[\nu:\hat C\to\Sigma,{\boldsymbol p}]$
to $B^{\,eg}_C$ gives rise to a bijection
between the set of elements $[\nu':\hat C'\to
\Sigma,{\boldsymbol p}']$ such that $\nu'
({\boldsymbol p}') = {\boldsymbol w}'$
on one side
and
$\Lambda({\boldsymbol w}')\cap B^{\,eg}_C$ on the other side.
\end{enumerate}
\end{lemma}

{\bf Proof.} The existence of the required set ${\mathcal C}$ immediately follows from the
relation \begin{equation}
h^0(C,{\mathcal J}^{cond}_C\otimes
{\mathcal O}_\Sigma(D))=n\ ,\label{eelra7}\end{equation} since imposing one by one $n$ generic point
constraints, we reduce $h^0$ to zero. To prove (\ref{eelra7}) we use the fact that
${\mathcal J}^{cond}_C=\nu_*{\mathcal O}_{\hat C}(-
\Delta)$, where $\Delta\subset\hat C$ is the so-called double-point divisor,
whose degree is $\deg\Delta=2\sum_{z\in\Sing(C)}\delta(C,z)$
(see, e.g., \cite[Section 2.4]{Caporaso_Harris:1998} or
\cite[Section 4.2.4]{Dolgachev:2013}). Thus,
$$\deg({\mathcal J}^{cond}_C\otimes
{\mathcal O}_\Sigma(D))=D^2-2\sum_{z\in\Sing(C)}\delta(C,z)=-DK_\Sigma+2g-2>2g-2\ ,$$ and hence
$$h^1(\hat C,{\mathcal J}^{cond}_C\otimes
{\mathcal O}_\Sigma(D))=0\quad\text{and}\quad h^0(\hat C,{\mathcal J}^{cond}_C\otimes
{\mathcal O}_\Sigma(D))=-DK_\Sigma+2g-2-g+1=n\ .$$

The dimension and the transversality in statement (1) mean
that the pull-back of $T_0B^{\,eg}_C$ to $|D|$
intersects
$|D|_{\boldsymbol w}$  transversally
and only at one point, which, in view of the
the identification of $T_0B^{\,eg}_C$ with $\prod_{z \in
\Sing(C)}J^{cond}_{C,z}/J_{C,z}$, reduces to (\ref{eeg}),
since ${\mathcal J}^{cond}_C$ can equivalently be regarded as
the ideal sheaf of
the zero-dimensional subscheme of $C$,
defined at all singular points $z\in\Sing(C)$ by the local conductor
ideals $J^{cond}_{C,z}$.

(2) The second statement of Lemma immediately follows from the
first one.
\proofend

\section{Appendix B: CH-configurations of points on real uninodal del Pezzo surfaces}\label{sec-ab}

Let $\Sigma$ be a uninodal del Pezzo surface of degree $\ge2$. Pick an effective divisor
class $D\in\Pic(\Sigma)$ represented by a curve not containing $E_\Sigma$ as a component,
and choose integer $g\ge0$ and two vectors $\alpha,\beta\in{\mathbb Z}^\infty_+$ such that
$I(\alpha+\beta)=DE_\Sigma$. Fix a sequence ${\boldsymbol w}$ of $\|\alpha\|$ distinct points
in general position on $E_\Sigma$ and a positive function $T:{\boldsymbol w}\to{\mathbb Z}$ such that
$|T^{-1}(i)|=\alpha_i$, $i\ge1$. Denote by
${\mathcal M}'_{g,\|\alpha\|}(\Sigma,D,\alpha,\beta,{\boldsymbol w},T)$ the space of elements
$[\nu:\hat C\to\Sigma,{\boldsymbol p}]\in{\mathcal M}'_{g,\|\alpha\|}(\Sigma,D,{\boldsymbol w})$ such that
\begin{itemize}\item $\nu^*({\boldsymbol w})=\sum_{p\in{\boldsymbol p}}T(\nu(p))\cdot p$,
\item $\nu^*(E_\Sigma\setminus{\boldsymbol w})=\sum_{q\in\hat C\setminus{\boldsymbol p}}k_q\cdot q$, where
the number of the coefficients $k_q$ equal to $i$ is $\beta_i$ for all $i\ge1$.
\end{itemize}

\begin{lemma}\label{lem-sm}
If ${\mathcal M}'_{g,\|\alpha\|}(\Sigma,D,\alpha,\beta,{\boldsymbol w},T)\ne\emptyset$, then
\begin{enumerate}\item[(i)] $\idim{\mathcal M}'_{g,\|\alpha\|}(\Sigma,D,\alpha,\beta,{\boldsymbol w},T)
\le n=-D(K_\Sigma+E_\Sigma)+g+\|\beta\|-1$,
\item[(ii)] a general element $[\nu:\hat C\to\Sigma,{\boldsymbol p}]$ of any
component of
${\mathcal M}'_{g,\|\alpha\|}(\Sigma,D,\alpha,\beta,{\boldsymbol w},T)$ of intersection dimension $n$ is an immersion,
and the curve $C=\nu(\hat C)$ is nonsingular along $E_\Sigma$.
\end{enumerate}
\end{lemma}

{\bf Proof.}
The statement
follows from \cite[Proposition 2.1]{Shoval_Shustin:2013}.
\proofend

Now let $\Sigma$ be a real uninodal del Pezzo surface
with a real $(-2)$-curve $E_\Sigma$
such that ${\mathbb R}E_\Sigma\ne\emptyset$. Pick an effective divisor class
$D_0\in\Pic(\Sigma)$, represented by a real curve not containing $E_\Sigma$ as component,
and such that $N=\dim|D_0|>0$. Denote by
$\Prec(D_0)$ the set of real effective divisor classes
$D\in\Pic(\Sigma)$, represented by real curves not containing $E_\Sigma$ as component,
and such that $D_0\ge D$. Notice that $\dim|D|\le N$.

Let $z_1,...,z_N$ be a sequence of $N$ distinct points in general position on ${\mathbb R}E_\Sigma$, and let
$z_i(t)$, $t\in[0,1]$, be a smooth path in ${\mathbb R}\Sigma$ transversal to
${\mathbb R}E_\Sigma$ at $z_i(0)=z_i$, $i=1,...,N$. We shall construct a sequence of
points $w_i=z_i(t_i)$, $0<t_i<1$, $i=1,...,N$, called a $D_0$-CH-configuration ({\it cf.}
\cite[Section 3.5.2]{Itenberg_Kharlamov_Shustin:2013}). We perform the construction inductively on $k=1,...,N$. Assume
that we have defined $t_i$, $i<k$, and then construct $t_k$ in the following procedure.
Given any data $D\in\Prec(D_0)$, $g\ge0$, $\alpha,\beta\in{\mathbb Z}^\infty_+$ such that
$I(\alpha+\beta)=DE_\Sigma$ and $1\le n=-D(K_\Sigma+E_\Sigma)+g+\|\beta\|-1\le k$, and given any
subsets $J_1\subset\{1,...,k-1\}$, $J_2\subset\{k+1,...,N\}$ such that
$|J_1|=n-1$, $|J_2|=\|\alpha\|$, we impose the following condition:

\noindent
for $t\in(0,t_k]$, the sets
\begin{equation}
\{[\nu:\hat C\to\Sigma,{\boldsymbol p}]\in{\mathcal M}'_{g,\|\alpha\|}
(\Sigma,D,\alpha,\beta,\{z_i\}_{i\in J_2},T)\ |\ w_i\in\nu(\hat C),\ i\in J_1,\ z_k(t)\in\nu(\hat C)\}\ ,\label{eelra100}\end{equation}
are finite of a capacity independent of $t$, and all their elements are presented by immersions.
The existence of such $t_k\in(0,1)$ follows from the fact that there are only finitely many
tuples $(D,g,\alpha,\beta,J_1,J_2)$, for which the sets (\ref{eelra100}), considered for
arbitrary $t\in(0,1)$, are
nonempty.

\begin{lemma}\label{lem-ch}
In the above notations, let ${\boldsymbol w}$ be a $D_0$-CH-configuration.
Suppose that $D\in\Prec(D_0)$, $g\ge0$, $\alpha,\beta\in
{\mathbb Z}^\infty_+$ satisfy $I(\alpha+\beta)=DE_\Sigma$ and $0\le n(D,
g,\beta)=-D(K_\Sigma+E_\Sigma)+g+\|\beta\|-1
\le N$. Then, for any disjoint sets $J_1,J_2\subset\{1,...,N\}$ such that
$|J_1|=n(D,g,\beta)$, $|J_2|=\|\alpha\|$, $\max J_1<\min J_2$, and for any
real element of the set
\begin{equation}\{[\nu:\hat C\to\Sigma,{\boldsymbol p}]\in{\mathcal M}'_{g,\|\alpha\|}
(\Sigma,D,\alpha,\beta,\{z_i\}_{i\in J_2},T)\ |\ w_i\in\nu(\hat C),\ i\in J_1\}
\label{eset}\end{equation}
the divisor $\nu^*(E_\Sigma)\subset\hat C$ is supported at only real points.
\end{lemma}

{\bf Proof.}
We use induction on $n=n(D,g,\beta)$. The case $n=0$ necessarily yields $g=0$ (see
\cite[Proposition 2.5]{Shoval_Shustin:2013}), and the desired statement follows then from
\cite[Lemma 3]{Itenberg_Kharlamov_Shustin:2013}. If $n>0$, we pick $k=\max J_1$ and consider degenerations
of a real element of the set (\ref{eset}) in the family
$[\nu_t:\hat C_t\to\Sigma,{\boldsymbol p}_t]$, $t\in(0,t_k]$, corresponding to the specialization of the point $w_k$ to $z_k\in E_\Sigma$ along the arc $L_k$.
By \cite[Proposition 2.6]{Shoval_Shustin:2013}, the limit of this family is
\begin{itemize}\item either $[\nu_0:\hat C_0\to\Sigma,{\boldsymbol p}_0]\in
{\mathcal M}'_{g,\|\alpha\|+1}(\Sigma,D,\alpha+e_m,\beta-e_m,\{z_i\}_{i\in J_2\cup\{k\}},T_0)$,
where $T_0\big|_{J_2}=T$, $T_0(z_k)=m$, which geometrically means that one of the nonfixed
intersection points of $\nu(\hat C)$ with $E_\Sigma$ of multiplicity $m$ becomes fixed
at the position $z_i$; the limit element satisfies
$n(D,g,\beta-e_m)=n-1$; hence by the induction assumption all intersection points
$\nu_0(\hat C_0)\cap E_\Sigma$ are real, and so are the points of $\nu(\hat C)\cap E_\Sigma)$;
\item or $[\nu_0:\hat C_0\to\Sigma,{\boldsymbol p}_0]$ such that $\hat C_0$ splits into components
$E_0$, $\hat C_1,...,\hat C_m$ so that $\nu_0:E_0\to E$ is an isomorphism, the elements
$[\nu_0:\hat C_j\to\Sigma,{\boldsymbol p}_j]\in{\mathcal M}'_{g_j,\|\alpha^{(j)}\|}(\Sigma,D_j,\alpha^{(j)},\beta^{(j)},T_j)$ satisfy
$n(D_j,g_j,\beta_j)<n$, $\nu_0({\boldsymbol p}_j)\subset\{z_i\ |\ i\in J_2\}$ for all $j=1,...,m$,
and, moreover, the divisor $\nu_t^*(E_\Sigma)$ is supported at the (real) points $\nu_t^{-1}(z_i)$, $1\le
i\le N$, and in a slightly deformed proper subset of the set
$\nu_0^{-1}(E_\Sigma\setminus\{z_1,...,z_N\})$; thus, if
$\nu_0:\hat C_j\to\Sigma$ is real, then by the induction assumption,
$\nu^{-1}_0(E_\Sigma)\cap\hat C_j$ is real, and hence a small smooth deformation of any of its proper subsets
is real too; if $\nu_0:\hat C_j\to\Sigma$ is not real (particularly, its complex conjugate
must be present in the splitting as well), then we have $n=0$ and $g=0$, which by
\cite[Lemma 3]{Itenberg_Kharlamov_Shustin:2013} implies that $E_\Sigma\cap\nu_0(\hat C_j)$ is just one point $z$; furthermore, in the deformation,
this node smooths out, and the deformed curve does not intersect $E_\Sigma$
in a neighborhood of $z$.
\end{itemize}
\proofend

%

{\it Address}: School of Mathematical Sciences, Tel Aviv
University, Ramat Aviv, 69978 Tel Aviv, Israel.
{\it E-mail}: {\tt shustin@post.tau.ac.il}

\end{document}